\documentclass[10pt]{amsart}

\usepackage{amsmath}
\usepackage{pxfonts}
\usepackage{amssymb}
\usepackage{amsfonts}

\newtheorem{theorem}{Theorem} [section]
\newtheorem{lemma} [theorem] {Lemma}

\newtheorem{corollary} [theorem] {Corollary}
\newtheorem{definition} [theorem] {Definition}
\newtheorem{remark} [theorem] {Remark}
\renewenvironment{proof} { \emph{Proof.} } { \rule{2mm}{2mm} \\}
\newenvironment{proof2} { \emph{Proof of Lemma \ref{dirichletexist}.} }{ \rule{2mm}{2mm} \\}
\newenvironment{proof3} { \emph{Proof of Theorem \ref{genT(B)}.} }{ \rule{2mm}{2mm} \\}
\newenvironment{proof4} { \emph{Proof of Lemma \ref{neumannexist}.} }{ \rule{2mm}{2mm} \\}
\newenvironment{proof5} { \emph{Proof of Theorem \ref{edge}.} }{ \rule{2mm}{2mm} \\}

\newcommand{\R}{{\bf R}}
\newcommand{\N}{{\bf N}}

\renewcommand{\b}{{\bf b}}

\newcommand{\n}{\nu}
\renewcommand{\t}{\tau}
\newcommand{\x}{X}
\newcommand{\y}{Y}
\newcommand{\z}{Z}
\newcommand{\0}{0}

\newcommand{\q}{Q}

\newcommand{\M}{{\mathcal{M}}}
\renewcommand{\S}{{\mathcal{S}}}

\newcommand{\esssup}{\mbox{ess sup}}

\newcommand{\tr}{\mbox{Tr}}

\newcommand{\dist}[2]{\mbox{dist}(#1,#2)}
\renewcommand{\div}{\mbox{div}}
\renewcommand{\leq}{\leqslant}
\renewcommand{\geq}{\geqslant}
\renewcommand{\iint}{\int}
\renewcommand{\Im}{\mbox{Im}}

\author[C.E. Kenig]{Carlos E. Kenig*}
\author[D.J. Rule]{David J. Rule}
\address{Department of Mathematics \\ The University of Chicago \\ 5734 S. University Avenue \\ Chicago \\ Illinois 60637 \\ USA}
\thanks{*Supported in part by the NSF}
\title[The Regularity and Neumann Problem]{The Regularity and Neumann Problem
for Non-symmetric Elliptic Operators}
\date{October 23, 2006 \newline
}

\begin{document}

\maketitle

\section{Introduction}

We will consider the Dirichlet problem
\begin{equation} \label{dirichlet}
\left\{ \begin{array}{ll} Lu = 0, & \mbox{in} \hspace{0.2cm}\Omega \\
u=f_0, & \mbox{on} \hspace{0.2cm} \partial \Omega
\end{array} \right. 
\end{equation}
with boundary data $f_0$ and the Neumann problem
\begin{equation} \label{neumann}
\left\{ \begin{array}{ll} Lu = 0, & \mbox{in} \hspace{0.2cm}\Omega \\
\n \cdot A\nabla u=g_0, & \mbox{on} \hspace{0.2cm} \partial \Omega
\end{array} \right. 
\end{equation}
with boundary data $g_0$. Here $\n$ is the outward unit normal vector to
$\partial\Omega$ and $L = \div \, A\nabla \cdot$  is an elliptic operator in
divergence form with coefficient matrix $A = (a_{ij})_{ij}$. The matrix $A$ is
assumed to have real-valued
bounded measurable entries ($\max_{i,j} \|a_{ij}
\|_{L^\infty(\Omega)} = \Lambda < \infty$) and satisfy the uniform ellipticity condition
\begin{equation} \label{ellipticity}
\lambda |\xi|^2 \leq \xi \cdot A\xi
\end{equation}
for some $\lambda > 0$ and all $\xi \in \R^2$, but $A$ is not necessarily symmetric. The domain $\Omega = \{ X = (x,t)
\in \R^2 \, | \, \phi(x) < t \}$ is the domain above the graph of a Lipschitz
function $\phi$. In the sequel we will denote by $\t$ the tangent $(1, \phi')/(1
+ (\phi')^2)^\frac{1}{2}$ to
$\partial\Omega$ and $\partial_\t = \t \cdot \nabla$ the derivative along the
boundary. The conormal derivative will be $\n \cdot A\nabla$.

Let us begin by fixing some notation. The spaces $C^\infty$ and $C^\infty_0$ are the spaces of smooth functions and smooth functions with compact support, respectively. For clarity, we will often indicate the domain of such functions, for example $C^\infty_0(\R)$. For $p \in [1,\infty]$ and $E$ a subset of
either $\R$ or $\R^2$ the space $L^p(E)$
is the set of measurable functions
$f \colon E \to \R$ such that $\| f \|_{L^p(E)}$ is finite, where
\[
\| f \|_{L^p(E)} = \left\{ \begin{array}{ll}\left( \int_{E} |f|^p
\right)^\frac{1}{p} , & \mbox{if} \hspace{0.2cm} p \in
[1,\infty) \\
\esssup_{E} |f|, & \mbox{if} \hspace{0.2cm} p=\infty
\end{array} \right.
\]
defined with respect to the appropriate Lebesgue measure. We define $L^p(E, \M)$ similarly for functions $F$ with values in $\M$, the set of
real-valued $2 \times 2$ matrices, by replacing the absolute value with $|F| :=
\sup_{ij}|f_{ij}|$ where $F = (f_{ij})_{ij}$.
The set $W^{1,p}(E)$ is the familiar Sobolev space, consisting of functions in
$L^p(E)$ whose first-order derivatives (in the sense of distributions) also
belong to $L^p(E)$ and $W^{1,p}_{\mbox{\tiny loc}}(E)$ the set consisting of those
functions in $W^{1,p}(E')$ for every compact subset $E'$ of $E$. We will consider
the
\emph{non-tangential approach
regions}
\[
\Gamma(\q) = \{ \x \in \Omega \,|\, |\x - \q| \leq (1 + a)
\dist{\x}{\partial\Omega} \}
\]
for each $\q \in \partial\Omega$ ($a > 0$ fixed). Here $\dist{\x}{\partial\Omega} = \inf_{\q \in
\partial\Omega} |\q - \x|$. Recall the \emph{non-tangential
maximal function} for a function $u$ on $\Omega$ is a function $N(u) \colon
\partial\Omega \to \R$ given by
\[
N(u)(\q) = \sup_{\Gamma(\q)} |u|
\]
and the related version
\[
\widetilde{N}(u)(\q) = \sup_{\x \in \Gamma(\q)} \left(
\frac{1}{|B_{\delta(\x)/2}(\x)|}\int_{B_{\delta(\x)/2}(\x)} |u|^2 \right)^\frac{1}{2}.
\]
Here $|E|$ is the Lebesgue measure of a set $E$ and $B_r = B_r(\x) = \{ \y \, | \,
|\x - \y| < r \}$ is the ball
centred at $\x$ of radius $r$. The notation $\delta(\x)$ is an abbreviation
for $\dist{\x}{\partial\Omega}$.

To solve (\ref{dirichlet}) and (\ref{neumann}) we will use the following weak formulations. It will be useful to introduce the space 
$\widetilde{W}^{1,2}(\Omega)$ of functions 
$f$ for which the norm
\begin{equation} \label{norm1}
\|f\|_{\widetilde{W}^{1,2}(\Omega)}:= \left( \iint_{\Omega} |f(\x)|^2 \, \frac{d\x}{(1 + |\x|^2)} + \iint_{\Omega} |\nabla f(\x)|^2 \, d\x \right)^\frac{1}{2}
\end{equation}
is finite and the space 
$\widetilde{L}^2(\partial\Omega)$ of functions 
$f$ for which the norm
\[
\|f\|_{\widetilde{L}^{2}(\partial\Omega)}:= \left(\iint_{\partial\Omega} |f(\x)|^2 \, \frac{d\sigma(\x)}{(1 + |\x|)} \right)^\frac{1}{2}
\]
is finite, where 
$d\sigma$ is Lebesgue surface measure on 
$\partial\Omega$. Denote by 
$\tr \colon \widetilde{W}^{1,2}(\Omega) \to \widetilde{L}^2(\partial\Omega)$ the trace operator initially defined on smooth functions 
$\varphi$ as
$\tr(\varphi) = \varphi|_{\partial\Omega}$ and 
$\widetilde{W}^{1,2}_0(\Omega)$ the set of all
$f \in \widetilde{W}^{1,2}(\Omega)$ with 
$\tr(f) = 0$. (One may readily check that $\tr$ has a continuous extention to $\widetilde{W}^{1,2}(\Omega)$, just as one shows $\tr \colon W^{1,2}(\Omega) \to L^2(\partial\Omega)$ is bounded (see, for example, \cite[p258]{E}).) 
\begin{lemma} \label{dirichletexist}
Given $f_0 \in W^{1,2}(\partial\Omega) \cap L^{\frac{7}{6}}(\partial\Omega) \cap L^{\frac{17}{6}}(\partial\Omega)$, there exists a unique 
$u \in \widetilde{W}^{1,2}(\Omega)$ such that 
$\tr(u) = f_0$ and
\begin{equation} \label{bilinear1}
\iint_{\Omega} A\nabla u \cdot \nabla \varphi = 0
\end{equation}
for all 
$\varphi \in \widetilde{W}^{1,2}_0(\Omega)$. Moreover, there exists a constant 
$C$, depending only on 
$\lambda$, 
$\Lambda$ and 
$\Omega$, such that 
\begin{equation} \label{space2}
\|u\|_{\widetilde{W}^{1,2}(\Omega)} \leq C\left(\|f_0\|_{W^{1,2}(\partial\Omega)} + \|f_0\|_{L^{7/6}(\partial\Omega)} + \|f_0\|_{L^{17/6}(\partial\Omega)}\right).
\end{equation}
\end{lemma}

We call the
$u$ from Lemma \ref{dirichletexist} the \emph{solution to (\ref{dirichlet}) with data 
$f_0$}. Let $H^1(\partial\Omega)$ denote the classical Hardy space $H^1(\R)$ (see, for example, \cite[p87]{S}) projected onto $\partial\Omega$, that is,
\[
H^1(\partial\Omega) = \{ f \colon \partial\Omega \to \R \, | \, x \mapsto f(x,\phi(x))(1 + \phi'(x)^2)^\frac{1}{2} \in H^1(\R) \}
\]
equiped with the obvious norm. We also introduce the homogeneous Sobolev space
\[
\dot{W}^{1,2}(\Omega) = \{ u \colon \Omega \to \R \, | \, \|\nabla u\|_{L^2(\Omega)} < \infty \},
\]
with norm $\| \cdot \|_{\dot{W}^{1,2}(\Omega)} = \|\nabla\cdot\|_{L^2(\Omega)}$.
\begin{lemma} \label{neumannexist}
Given 
$g_0 \in H^1(\partial\Omega)$, there exists a unique 
$u \in \dot{W}^{1,2}(\Omega)$ such that
\[
\iint_{\Omega} A\nabla u \cdot \nabla\varphi = \int_{\partial\Omega} g_0 \tr(\varphi) \, d\sigma
\]
for all 
$\varphi \in \dot{W}^{1,2}(\Omega)$. Moreover, there exists a constant 
$C$, depending only on 
$\lambda$, 
$\Lambda$ and 
$\Omega$, such that
\begin{equation} \label{space3}
\|u\|_{\dot{W}^{1,2}(\Omega)} \leq C\|g_0\|_{H^1(\partial\Omega)}.
\end{equation}
\end{lemma}

We call the
$u$ from Lemma \ref{neumannexist} the \emph{solution to (\ref{neumann}) with data 
$g_0$}. We will postpone the proofs of these lemmata until Section \ref{prelim}. Since, in general, we cannot assign a meaning to 
$\n \cdot A\nabla u|_{\partial\Omega}$, we should be careful as to exactly how we interpret the statement in (\ref{neumann}) that $\n \cdot A\nabla u = g_0$ on 
$\partial\Omega$. It is well known that existence of the estimates in the following definition enable a certain non-tangential convergence to the boundary data to be established (see, for example, \cite{KP}).

\begin{definition}
(i) We say that the Dirichlet problem holds for $p$, or $(D)_p^A = (D)_p$ holds, if for any $u$ solving (\ref{dirichlet}) with
boundary data $f_0 \in L^p(\partial\Omega) \cap W^{1,2}(\partial\Omega) \cap L^{\frac{7}{6}}(\partial\Omega) \cap L^{\frac{17}{6}}(\partial\Omega)$ we have
\[
\|N(u)\|_{L^p(\partial\Omega)} \leq C(p) \|f_0\|_{L^p(\partial\Omega)}.
\]
(ii) We say that the Neumann problem holds for $p$, or $(N)_p^A = (N)_p$ holds,
if for any $u$ solving (\ref{neumann}) with
boundary data $g_0 \in L^p(\partial\Omega) \cap H^1(\partial\Omega)$ we have
\[
\|\widetilde{N}(\nabla u)\|_{L^p(\partial\Omega)} \leq C(p) \|g_0\|_{L^p(\partial\Omega)}.
\]
(iii) We say that the regularity problem holds for $p$, or $(R)_p^A = (R)_p$ holds, if for any $u$ solving (\ref{dirichlet}) with
boundary data $f_0 \in W^{1,p}(\partial\Omega) \cap W^{1,2}(\partial\Omega) \cap L^{\frac{7}{6}}(\partial\Omega) \cap L^{\frac{17}{6}}(\partial\Omega)$ we have
\[
\|\widetilde{N}(\nabla u)\|_{L^p(\partial\Omega)} \leq C(p) \|\partial_\t f_0\|_{L^p(\partial\Omega)}.
\]
In each case, the constant $C(p) > 0$ must depend only on $\lambda$, $\Lambda$,
$\Omega$ and
$p$.
\end{definition}

The main result of this paper is the following.

\begin{theorem} \label{main}
Let $L = \div \, A \nabla$ be an elliptic operator as defined above with
coefficient matrix $A = A(x)$ independent of the $t$-variable in the domain
$\Omega = \{ X = (x,t) \in \R^2 \, | \, \phi(x) < t \}$ above the graph of a
Lipschitz function $\phi$. Then $(N)_p$ and $(R)_p$ hold for some (possibly small) $p >
1$.
\end{theorem}

In 2000, Kenig, Koch, Pipher and Toro~\cite{KKPT} showed that the Dirichlet
problem for an elliptic operator $L = \div \, A\nabla$ in
$\Omega \subset \R^2$ holds for some $p \geq 2$, where $A$ has coefficients
independent of the  $t$-direction. The advance here was that the $2 \times 2$
matrix $A$ did not need to be assumed symmetric. The exponent $p$ could not in
general be specified, as they showed with an example. This example can be used to show the same is true in our case: For any given $p$ there exist operators $L$ as in Theorem \ref{main} for which $(N)_p$ and $(R)_p$ do not hold. (This example is worked out in the Appendix.) Jerison and Kenig showed in
\cite{JK1} that when $A$ is also symmetric $(D)_2$ holds.
In \cite{KP} Kenig and Pipher show that if $(R)_p^A$ holds in $\Omega$
for some $p$ then $(D)_{p'}^{A^t}$ holds, where $\frac{1}{p} +
\frac{1}{p'} = 1$ and $A^t$ is the transpose of $A$. Here we will prove Theorem \ref{main} by showing a reverse
implication: if $(D)_{p'}^{A^t}$ holds and $A$ is independent of the
$t$-direction, then $(N)_p^{\widetilde{A}}$ and $(R)_p^A$ hold, where
$\widetilde{A} = A^t/\det(A)$. Then, as a simple corollary, the main result in \cite{KKPT}
passes over to the Neumann and regularity problems to give Theorem \ref{main}.

The layout of this paper is as follows. In Section \ref{prelim} we fix some more
notation and lay out some preliminaries. In Section \ref{reduct} we begin by
observing that
a change of variables used in \cite{KKPT} can reduce our problem to
coefficient matrices $A$ which are upper trianglular. We then go on to show that
the solvability of the Dirichlet problem allows us to reduce Theorem \ref{main}
to the boundedness of certain layer potentials. In Section \ref{boundedness},
this boundedness is proved in the special case that the boundary of our domain
has a small Lipschitz constant. In Section \ref{david} we use David's method
\cite{D} to show that the layer potentials must also be bounded on domains
$\Omega$ when $\phi$ is an arbitrary Lipschitz function. As an appendix, we show the example in \cite{KKPT} serves the same purpose for Theorem \ref{main} as it does for their Theorem 3.1.

As is common practice the letter $C$ will denote a constant whose value may
change from line to line, but which can be fixed depending only on the constants
$\lambda$, $\Lambda$ and the Lipschitz constant of $\phi$. When a subscript is
added the value will not vary from line to line, and when the constant depends
on other parameters or when we wish to record the form of the dependency this
will be noted explicity (for example, we wrote above $C = C(p)$).

We would like to thank Steve Hofmann for helpful conversations and for making the related
work by him and coauthors \cite{AAAHK} available to us.

\section{Preliminaries}

\label{prelim}

The \emph{square function} $S(f) \colon \partial\Omega \to \R$ of $f \colon
\Omega \to \R$ is given by
\[
S(f)(\q) = \left\{ \int_{\Gamma(\q)} |\nabla f|^2
\right\}^\frac{1}{2}
\]
and the \emph{Hardy-Littlewood maximal operator} $M \colon L^p(\partial\Omega)
\to L^p(\partial\Omega)$ ($1 < p \leq \infty$) is given by
\[
M(f)(\q) = \sup_{r > 0} \frac{1}{\sigma(\Delta_r(\q))} \, \int_{\Delta_r(\q)}
|f| \, d\sigma
\]
where $\Delta_r(\q) = B_r(\q) \cap \partial\Omega$ and $d\sigma$ is the Lebesgue
surface measure of $\partial\Omega$. We also define $M$ acting on
functions $F \colon \R \to \M$ analogously, replacing the absolute value with the corresponding
matrix version and taking the supremum over intervals in $\R$.
A function $K \colon \R^2 \to \M$ is said
to be a \emph{Calder\'on-Zygmund kernel} if there exists constants $C > 0$ and $\alpha
\in (0,1]$ such that
\begin{equation} \label{kerest1}
|K(x,y)| \leq \frac{C}{|x-y|}
\end{equation}
for all $x,y \in \R$, $x \neq y$,
\begin{equation} \label{kerest2}
|K(x,y) - K(x',y)| \leq \frac{C|x-x'|^\alpha}{(|x-y| + |x'-y|)^{1 + \alpha}}
\end{equation} 
when $|x-x'| \leq \frac{1}{2}\max\{|x-y|,|x'-y|\}$, and
\begin{equation} \label{kerest3}
|K(x,y) - K(x,y')| \leq \frac{C|y-y'|^\alpha}{(|x-y| + |x-y'|)^{1 + \alpha}}
\end{equation}
when $|y-y'| \leq \frac{1}{2}\max\{|x-y|,|x-y'|\}$.
Recall, for $F = (f_{ij})_{ij}$, $|F| := \sup_{ij}|f_{ij}|$.
For Banach spaces $\mathcal{X}$ and $\mathcal{Y}$ dense in $L^2(\R, \M)$, a
continuous linear operator $T \colon \mathcal{X} \to \mathcal{Y'}$ is said to be
a \emph{singular integral operator} associated to the Calder\'on-Zygmund kernel $K$ if for each $F
\in \mathcal{X}$ and $G \in \mathcal{Y}$ with disjoint support, we have the
representation
\begin{equation} \label{singrep}
\langle G, T(F) \rangle = \int_{\R^2} G(x)^tK(x,y)F(y) \, dydx,
\end{equation}
where $\langle \cdot , \cdot \rangle$ is the standard inner product on
$L^2(\R,\M)$ and $G(x)^t$
is the transpose of $G(x)$. The \emph{maximal singular integral operator}
$T^*$ associated to $K$ is defined by first setting
\[
T^{(\delta)}(F)(x) = \int_{|y-x| \geq \delta} K(x,y)F(y) \, dy
\]
and then $T^*(F) = \sup_{\delta > 0} |T^{(\delta)}(F)|$. The following lemmata are well known for any elliptic operator as defined above
(see, for example, \cite{K} and \cite[Thm 8.24]{GT}).

\begin{lemma} \label{cacciopoli}
Let $u \in H^{1,2}(B_{2r})$ satisfy $Lu = 0$, then there exists a constant $C > 0$ which depends only on
$\lambda$ and $\Lambda$ such that
\[
\fint_{B_r} |\nabla u(\z)|^2 \, d\z \leq \frac{C}{r^2} \, \fint_{B_{2r}} |u(\z)|^2 \,
d\z.
\]
\end{lemma}

\begin{lemma} \label{sup}
For $u \in H^{1,2}(B_{5r/4})$ such that $Lu = 0$ and $p \geq 1$, there exists a constant
$C(p) > 0$ which depends only on
$\lambda$, $\Lambda$ and $p$ such that
\[
\sup_{B_{r/2}} |u| \leq C(p) \left( \fint_{B_r} |u|^p \right)^\frac{1}{p}
\]
\end{lemma}

\begin{lemma} \label{reversehoelder}
For $u \in H^{1,2}(B_{3r/2})$ such that $Lu = 0$, there exist constants
$C > 0$ and $p > 2$ which depend only on
$\lambda$ and $\Lambda$ such that
\[
\left( \fint_{B_{r/2}} |\nabla u|^p \right)^\frac{1}{p} \leq C \left(
\fint_{B_r} |\nabla u|^2 \right)^\frac{1}{2}
\]
\end{lemma}

\begin{lemma} \label{Calpha}
Let $u \in H^{1,2}(B_r)$ satisfy $Lu = 0$, then there exists constants
$C > 0$ and $\alpha \in (0,1]$, depending only on $\lambda$ and $\Lambda$, such
that
\[
|u(\x) - u(\y)| \leq \frac{C}{r^{1+\alpha}} \|u\|_{L^2(B_r)}|\x - \y|^\alpha,
\]
for all $\x, \y \in B_{2r/3}$.
\end{lemma}

We will now prove Lemmata \ref{dirichletexist} and \ref{neumannexist}.

\begin{proof2}
Without loss of generality we can suppose $\Omega = \R^2_+ = \{(x,t) \, | \, t > 0\}$, since the transformation 
$\Phi \colon \Omega \to \R^2_+$ given by 
$\Phi(x,t) = (x,t-\phi(x))$ will yield the general case. Let 
$P_t \colon \R \to \R$ denote the Poisson kernel:
\[
P_t(x) := \frac{t}{\pi(|x|^2 + t^2)}.
\]
Let 
\[
w(x,t) = f_0 \ast P_t(x) := \int_\R f_0(x-y)P_t(y) \, dy
\]
be the harmonic extension of 
$f_0$ to 
$\R^2_+$. We will show that
\begin{equation} \label{space1}
\|w\|_{\widetilde{W}^{1,2}(\R^2_+)} \leq C(\|f_0\|_{W^{1,2}(\partial\R^2_+)} + \|f_0\|_{L^{7/6}(\partial\R^2_+)} + \|f_0\|_{L^{17/6}(\partial\R^2_+)}).
\end{equation}
Using H\"older's inequality
\[
\iint_{\R^2_+} |w(\x)|^2 \, \frac{d\x}{(1 + |\x|^2)} \leq C \|w\|_{L^{17/6}(\R^2_+)}^2,
\]
but, by Fubini's Theorem,
\[
\begin{aligned}
\|w\|_{L^{17/6}(\R^2_+)}^{\frac{17}{6}} & = \int_0^\infty \|f_0 \ast P_t\|_{L^{17/6}(\R)}^{\frac{17}{6}} \, dt \\
& =  \int_0^1 \|f_0 \ast P_t\|_{L^{17/6}(\R)}^{\frac{17}{6}} \, dt + \int_1^\infty \|f_0 \ast P_t\|_{L^{17/6}(\R)}^{\frac{17}{6}} \, dt.
\end{aligned}
\]
Now, using well-known properties of the Poisson kernel and Young's inequality (see, for example, \cite{S2}), we can bound this by
\[
\begin{aligned}
& \|M(f_0)\|_{L^{17/6}(\R)}^{\frac{17}{6}} + \|f_0\|_{L^{7/6}(\R)}^{\frac{17}{6}} \int_1^\infty \|P_t\|_{L^{119/59}(\R)}^{\frac{17}{6}} \, dt \\
& \leq C\|f_0\|_{L^{17/6}(\R)}^{\frac{17}{6}} + C \|f_0\|_{L^{7/6}(\R)}^{\frac{17}{6}},
\end{aligned}
\]
since 
$\|P_t\|_{L^{119/59}(\R)}^{\frac{17}{6}} = Ct^{-\frac{10}{7}}$. This controls the first integral in the norm 
$\|w\|_{\widetilde{W}^{1,2}(\R^2_+)}$. To control the second, first observe that
\[
\|P_t'\|_{L^1(\R)} = Ct^{-2},
\]
so, as before,
\[
\begin{aligned}
& \|\partial_1w\|_{L^2(\R^2_+)} \leq \int_0^1 \|f_0' \ast P_t\|_{L^2(\R)} \, dt + \int_1^\infty \|f_0 \ast P_t'\|_{L^2(\R)} \, dt \\
& \leq \|M(f_0')\|_{L^2(\R)} + \|f_0\|_{L^2(\R)} \int_1^\infty \|P_t'\|_{L^1(\R)} \, dt \leq C\|f_0\|_{H^{1,2}(\partial\R^2_+)}.
\end{aligned}
\]
We may control 
$\|\partial_2w\|_{L^2(\R^2_+)}$ similarly since (via the Cauchy-Riemann equations) we know 
$\partial_2w = -\partial_1\widetilde{w}$, where 
$\widetilde{w}(x,t) = H(f_0) \ast P_t(x)$ and 
$H \colon L^2(\partial\R^2_+) \to L^2(\partial\R^2_+)$ is the Hilbert transform. Thus (\ref{space1}) is proved.

Now we will show
\begin{equation} \label{bilinearform}
(\psi,\varphi) \mapsto \iint_{\R^2_+} A\nabla\psi \cdot \nabla \varphi
\end{equation}
defines a bounded coercive bilinear form on 
$\widetilde{W}_0^{1,2}(\R^2_+)$. Linearity and boundedness are clear, so to prove coercivity we first observe that 
$\psi \in \widetilde{W}_0^{1,2}(\R^2_+)$ can be extended to a function in 
$\widetilde{W}^{1,2}(\R^2)$ (the space with norm (\ref{norm1}) but integrating over 
$\R^2$), which we will also call 
$\psi$, simply by extending by zero. Then obviously
$\|\nabla\psi\|_{L^2(\R^2)} = \|\nabla\psi\|_{L^2(\R^2_+)}$.
Thus, using Poincar\'e's inequality,
\[
\begin{aligned}
\iint_{\R^2} |\psi(\x)|^2 \, \frac{d\x}{(1 + |\x|^2)} & \leq C\int_{B_4(\0)} |\psi|^2 + C\sum_{n=2}^\infty 2^{-2n}\int_{B_{2^n}(\0) \setminus B_{2^{n-1}}(\0)} |\psi|^2 \\
& \leq C\int_{B_4(\0)} |\nabla \psi|^2 + C\sum_{n=2}^\infty \int_{B_{2^n}(\0) \setminus B_{2^{n-1}}(\0)} |\nabla \psi|^2 \\
& \leq C \iint_{\R^2} |\nabla \psi|^2 = C \iint_{\R^2_+} |\nabla \psi|^2.
\end{aligned}
\]
Combining this with (\ref{ellipticity}) easily gives the coercivity.

We may now apply the Lax-Milgram Theorem to find a unique $v \in \widetilde{W}^{1,2}_0(\R^2_+)$ such that
\[
\iint_{\R^2_+} A\nabla v \cdot \nabla\varphi = \iint_{\R^2_+} -A\nabla w \cdot \nabla \varphi
\]
for all 
$\varphi \in \widetilde{W}^{1,2}_0(\R^2_+)$. Then 
$u = v + w \in \widetilde{W}^{1,2}(\R^2_+)$ is clearly the unique function in 
$\widetilde{W}^{1,2}(\R^2_+)$ such that (\ref{bilinear1}) holds and 
$\tr(u) = f_0$. The estimate (\ref{space2}) follows from (\ref{space1}) and the Lax-Milgram Theorem.
\end{proof2}

A function $f \colon \R \to \R$ is said to be of \emph{bounded mean
oscillation}, written $f \in \mbox{BMO}$, if
\[
\|f\|_{\mbox{\tiny BMO}} := \sup_I \fint_I |f(x) - \fint_I f| \, dx < \infty,
\]
where the supremum is taken over all intervals $I \in \R$.

\begin{proof4}
Once again, without loss of generality we may assume $\Omega = \R^2_+$. The bilinear form (\ref{bilinearform}) is clearly bounded and coercive on $\dot{W}^{1,2}(\R^2_+)$ so to apply the Lax-Milgram Theorem it remains to check that
\[
\varphi \mapsto \int_{\partial\Omega} g_0 \tr(\varphi) \, d\sigma
\]
is a bounded linear functional on $\dot{W}^{1,2}(\R^2_+)$ for a given $g_0 \in H^1(\partial\R^2_+)$.

It suffices to check $\tr \colon \dot{W}^{1,2}(\R^2_+) \to \mbox{BMO}$ is a bounded operator. This is easily done as follows. Let $I = B_R \cap \partial\R^2_+$ and $B^+_R = B_R \cap \R^2_+$, where $B_R$ is a ball with centre on $\partial\R^2_+$. Denote by $\varphi_{B^+_{2R}}$ the average $\fint_{B^+_{2R}} u$, where $B_{2R}$ is the ball concentric with $B_R$, but with radius $2R$. Fix $\xi \colon \overline{\R^2_+} \to \R$ to be a smooth cut-off function equal to one near $I$, supported in $B_{2R}$ and such that $|\nabla \xi| \leq C/R$. Then
\[
\begin{aligned}
\fint_I |\varphi - \varphi_{B^+_{2R}}| \, d\sigma & \leq \frac{1}{R} \int_{\partial\R^2_+} \xi|\varphi - \varphi_{B^+_{2R}}| = - \frac{1}{R}\iint_{B^+_{2R}} \partial_t(\xi|\varphi - \varphi_{B^+_{2R}}|) \\
& = - \frac{1}{R}\iint_{B^+_{2R}} (\partial_t\xi)|\varphi - \varphi_{B^+_{2R}}| - \frac{1}{R}\iint_{B^+_{2R}} \xi(\partial_t\varphi)\mbox{sgn}(\varphi - \varphi_{B^+_{2R}}) \\
& \leq \frac{C}{R^2}\iint_{B^+_{2R}} |\varphi - \varphi_{B^+_{2R}}| + \frac{1}{R}\iint_{B^+_{2R}} |\nabla \varphi| \\
& \leq C \|\nabla \varphi \|_{L^2(\R^2_+)} = C \|\varphi \|_{\dot{W}^{1,2}(\R^2_+)},
\end{aligned}
\]
where the last inequality follows from H\"older's and Poincar\'e's inequalities.
\end{proof4}

Kenig and Ni \cite{KN} provide the following definition and existence of a
fundamental solution.
\begin{definition} \label{fundamentaldef}
A function $\Gamma_{\x} \colon \R^2 \to \R$ is called a
fundamental solution for $L = \div \, A\nabla \cdot$ with pole at $\x$ if
\begin{enumerate}
\item[(i)] $\Gamma_{\x} \in W^{1,2}_{\mbox{\tiny loc}}(\R^2 \setminus \{\x\}) \cap
W^{1,p}_{\mbox{\tiny loc}}(\R^2)$ for all $p < 2$, and, for
every $\varphi
\in C^\infty_0(\R^2)$,
\[
\iint_{\R^2} A^t(\y)\nabla \Gamma_{\x}(\y) \cdot \nabla \varphi(\y) \, d\y =
-\varphi(\x).
\]
\item[(ii)] $|\Gamma_\x(\y)| = O(\log|\x - \y|)$ as $|\y| \to \infty$.
\end{enumerate}
\end{definition}

\begin{theorem} \label{fundamental}
For each $\x \in \R^2$ there exists a unique fundamental solution $\Gamma_\x$
for $L$ with pole at $\x$, and positive constants $C_1, C_2, R_1 < 1,
R_2 > 1$, which depend only on $\lambda$ and $\Lambda$, such that
\[
\begin{aligned}
& C_1\log(1/|\x - \y|) \leq -\Gamma_\x(\y) \leq C_2 \log(1/|\x - \y|)
\hspace{0.3cm} \mbox{for} \hspace{0.3cm} 
 |\x - \y| <
R_1, \hspace{0.3cm} \mbox{and} \\
& C_1\log(|\x - \y|) \leq \Gamma_\x(\y) \leq C_2 \log(|\x - \y|) \hspace{0.3cm} \mbox{for}
\hspace{0.3cm} |\x - \y| >
R_2.
\end{aligned}
\]
\end{theorem}

The notation we employ to denote differentiation is standard, but let us clarify
one point. When
the fundamental solution is differentiated $\nabla \Gamma_\x(\y)$ will denote
the gradient in the $\y$-variable and $\nabla_\x\Gamma_\x(\y)$ the gradient in
the $\x$-variable.

\begin{lemma} \label{symmetry} Fix $\x, \y \in \R^2$. Let $\Gamma_\x$ be the
fundamental solution for an elliptic operator
$L$ as defined above with pole at $\x$, and $\Gamma^t_\y$ the fundamental
solution to the adjoint operator $L^t = \div \, A^t \nabla$ with pole at $\y$.
Then we have that
\[
\Gamma_\x(\y) = \Gamma^t_\y(\x).
\]
\end{lemma}

\begin{proof}
We will prove the lemma for a modified coefficient matrix $A$ where $A = I$
outside the ball $B_R(\0)$. It is then easy using, Definition
\ref{fundamentaldef} and Theorem
\ref{fundamental}, to
show that passing to the limit $R \to \infty$ we obtain the identity for
general $A$.

The linear transformation $\z \mapsto C \z$ ($C > 0$), transforms
$\Gamma_\x$ and $\Gamma_\y$ into fundamental solutions of operators with the
same ellipticity
constants $\lambda$ and $\Lambda$. Thus, without loss of generality, we can assume $A
= I$ outside $B_{1/2}(\0)$ and that $\x, \y \in B_{1/2}(\0)$.
Green's second identity gives us
\[
\begin{aligned}
\Gamma^t_\y(\x) - \Gamma_\x(\y) & = \int_{\partial B_R(\0)} \Gamma^t_\y(\z)
\n(\z) \cdot A^t(\z) \nabla
\Gamma_\x(\z) \, d\sigma(\z) \\
& - \int_{\partial B_R(\0)} \Gamma_\x(\z) \n(\z)
\cdot A(\z) \nabla
\Gamma^t_\y(\z) \, d\sigma(\z),
\end{aligned}
\]
for any $R > 1$. But on the other hand, the uniqueness and the construction used
in \cite{KN} shows us that, for $\z \in B_1(\0)^c$,
\[
\Gamma_\x(\z) = g(\z) + w(\z)
\]
and
\begin{equation} \label{fundrep}
\Gamma^t_\y(\z) = h(\z) + v(\z),
\end{equation}
where $g$ and $h$ are Green's functions for $B_1(\x)^c$ and $B_1(\y)^c$
respectively, with poles at
infinity, and $w$ and $v$ are harmonic functions satisfying $|w(\z)| + |v(\z)|
\leq C|\z|^{-\alpha}$, so then $|\nabla w(\z)| + |\nabla v(\z)| \leq
C|\z|^{-(1+\alpha)}$ (see, for example, \cite[Thm 2.10]{GT}). Given our
assumptions on $A$ we have the explicit expressions $g(\z) = \log |\z - \x|$ and
$h(\z) = \log |\z - \y|$. Thus
\[
\begin{aligned}
|\nabla \Gamma_\x(\z) - \nabla \Gamma_\y^t(\z)| & \leq |\nabla g(\z) - \nabla h(\z)| + |\nabla w(\z)| +
|\nabla v(\z)| \\
& \leq C(R^{-2} + R^{-(1+\alpha)})
\end{aligned}
\]
and
\[
|\Gamma_\x(\z) - \log|\z|| + |\Gamma^t_\y(\z) - \log|\z|| \leq C(R^{-\alpha} +
\log|(R+1)/R|)
\]
on $\partial B_R(\0)$. And we may conclude, by the above and our assumptions on
$A$, that
\[
\begin{aligned}
|\Gamma^t_\y(\x) - \Gamma_\x(\y)| & \leq \int_{\partial B_R(\0)} |(\Gamma^t_\y(\z)
- (\log|\z|)) \n(\z) \cdot A(\z) \nabla
\Gamma_\x(\z)| \, d\sigma(\z) \\
& + \int_{\partial B_R(\0)} |(\log|\z|) \n(\z) \cdot (A^t(\z) \nabla
\Gamma_\x(\z) - A \nabla\Gamma^t_\y(\z))|  \, d\sigma(\z)\\
& + \int_{\partial B_R(\0)} |((\log|\z|) - \Gamma_\x(\z)) \n(\z)
\cdot A(\z) \nabla
\Gamma^t_\y(\z) \, d\sigma(\z) \\
& \leq C(R^{-\alpha} + \log|(R+1)/R| + (\log R)(R^{-\alpha} + R^{-1}))
\end{aligned}
\]
which tends to zero as $R \to \infty$.
\end{proof}




The following theorem was proved by Kenig, Koch, Pipher and Toro~\cite{KKPT}.
They prove that if the $L^2$-norms of the square function and non-tangential maximal
function of all solutions vanishing at a fixed point can be compared in all bounded Lipschitz domains
$\Omega_0$
contained in $\Omega$ with
constants depending only on the Lipschitz character of the domain, then
$(D)_p$ holds for
some $p$. The estimate
\begin{equation} \label{square}
\|S(u)\|_{L^p(\partial\Omega_0)} \simeq
\|N(u)\|_{L^p(\partial\Omega_0)}
\end{equation}
for each $p > 1$ and solutions $u$ which vanish at a fixed point in
$\Omega_0$ is demonstrated therein. Dahlberg \cite{D1} showed this estimate for
$p=2$ implies that solutions are $\varepsilon$-approximable and in \cite{KKPT}
the authors go on to show this $\varepsilon$-approximability implies that $(D)_p$ holds for some $p > 1$. We refer the reader to \cite{KKPT}
for the details and precise definitions, Remark (2.11) providing a brief overview.

\begin{theorem} \label{KKPT}
Let $L = \div \, A\nabla$ be an elliptic operator in a domain $\Omega = \{(x,t) \in \R^2
\, | \, \phi(x) < t \}$, where $A = A(x)$ is independent of the $t$-variable and
$\phi \colon \R \to \R$ is a Lipschitz function. Then
there exists a (possibly large) $p$ such that $(D)_p$ holds in $\Omega$, with bound depending only on $\lambda$, $\Lambda$, $p$ and the Lipschitz constant of $\phi$.
\end{theorem}

Given any solution $u$ to an elliptic equation $Lu = \div \, A
\nabla u = 0$, the vector $A \nabla u$ is divergence free and so can be written
as the curl of a vector. This amounts to finding a $\tilde{u}$ solving
\begin{equation} \label{conj}
\left(  \begin{array}{cc} 0 & 1 \\ -1 & 0 \end{array} \right)\nabla \tilde{u} = A\nabla
u.
\end{equation}
The function $\tilde{u}$ (defined up to a constant by (\ref{conj})) is called the \emph{conjugate} of $u$ and
$(u,\tilde{u})$ is a \emph{conjugate pair}. Observe firstly that $\widetilde{u}$
satisfies an elliptic equation with coefficient matrix $\widetilde{A} =
A^t/\det{A}$, and secondly that the conormal derivative
of $u$ is the tangential derivative of $\tilde{u}$ and vice versa.

In Sections \ref{reduct} and \ref{boundedness} we will work under the a priori assumptions that $A=I$ for large $x$, $A$ and $\phi$ are smooth functions, $\|\phi'\|_{L^\infty(\R)} \leq k$, $\phi' \equiv \alpha_0$ for large $x$ and $x \mapsto \phi(x) - \alpha_0x \in C^\infty_0(\R)$. Once our theorems have been proved under our a priori assumptions, it is a simple matter to obtain the general case. Note that, under our a priori assumptions, if $u$ solves (\ref{dirichlet}) with data $f_0 \in C^\infty_0(\partial\Omega)$, then $u \in C^\infty(\overline{\Omega})$, and $u(\x) = O(|\x|^{\delta-1})$ and $\nabla u(\x) = O(|\x|^{\delta-2})$ for all $\delta > 0$ as $|\x| \to \infty$. Moreover, if
\begin{equation} \label{conj2}
\widetilde{u}(\x) = - \int_{\gamma(\x)} \n(\z) \cdot A(\z)\nabla u(\z) \, dl(\z),
\end{equation}
where $\gamma(\x)$ is the line segment $\gamma(\x) = \{(y,s) \, | \, s \geq t, y=x, \x = (x,t) \}$, $\n$ is the unit normal to $\gamma(\x)$ and $dl$ is arc length, then $(u,\widetilde{u})$ is a conjugate pair and $\widetilde{u}(\x) = O(|\x|^{\delta-1})$ for all $\delta > 0$ as $|\x| \to \infty$. (See (\ref{gamtildedef}) and the lines following it for a similar argument.)

\begin{theorem} \label{edge}
Let $u$, $\Omega$ and $A$ be as in the previous paragraph. If $p' > 1$ is such that
$(D)^{A^t}_{p'}$
holds then there exists a
constant $C(p)$, depending only on $\lambda$,
$\Lambda$, $k$, $p$ and the $(D)^{A^t}_{p'}$ constant of $A^t$, such that
\begin{equation} \label{comp2}
\|\nabla u\|_{L^p(\partial\Omega)} \leq C(p) \|\partial_\t
f_0\|_{L^p(\partial\Omega)}.
\end{equation}
Also, if $u$ solves (\ref{neumann}) with coefficient matrix $A$ replaced by $\widetilde{A} = A^t/\det(A)$ and $u$ verifies the a priori assumptions $u(\x) = O(|\x|^{\delta-1})$ and $\nabla u(\x) = O(|\x|^{\delta-2})$ for all $\delta > 0$ as $|\x| \to \infty$, then there exists a constant $C(p)$, depending on the same quantities, such that
\begin{equation} \label{comp3}
\|\nabla u\|_{L^p(\partial\Omega)} \leq C(p) \|g_0\|_{L^p(\partial\Omega)}.
\end{equation}
As usual $\frac{1}{p} + \frac{1}{p'} = 1$.
\end{theorem}

We cannot expect these estimates
to hold for all $p > 1$ as the Appendix shows. The a priori assumptions for solutions to (\ref{neumann}) will hold if, for instance, $g_0 \in C^\infty_0(\partial\Omega) \cap H^1(\partial\Omega)$. Such $g_0$ are dense in $L^p(\partial\Omega)$, $1 < p < \infty$.

\begin{proof5}
We start out by showing that if $u(\x) = O(|\x|^{\delta-1})$ for all $\delta > 0$ as $|\x| \to \infty$, then, for $1 < p < \infty$,
\begin{equation} \label{square2}
\|S(u)\|_{L^p(\partial\Omega)} \simeq
\|N(u)\|_{L^p(\partial\Omega)}.
\end{equation}
In fact, let $\Omega_R = \Omega \cap B_R(0,\phi(0))$ so we have that $A_R = (0,R/2) \in \Omega_R$, $\dist{A_R}{\partial\Omega_R} \newline \simeq \dist{A_R}{\partial\Omega} \simeq R$. By (\ref{square}) we have
\[
\|S(u)\|_{L^p(\partial\Omega_R)} \simeq
\|N(u - u(A_R))\|_{L^p(\partial\Omega_R)},
\]
but since $u(\x) = O(|\x|^{\delta-1})$, $|\nabla u| = O(|\x|^{\delta-2})$ for all $\delta > 0$ we may choose $\delta$ so that $1 + (\delta - 1)p > 0$ so $|u(A_r)|^p|\partial\Omega_R| + \|S(u)\|_{L^p(\partial\Omega_R \setminus \partial\Omega)} + \|N(u)\|_{L^p(\partial\Omega_R \setminus \partial\Omega)} \to 0$ as $R \to \infty$ and (\ref{square2}) follows. Now let $f_0$, $u$ and $\widetilde{u}$ be as above but with $A$ replaced by $A^t$. Define the operator $H_{A^t}$ by $H_{A^t}(f_0) = \widetilde{u}|_{\partial\Omega}$.
This a bounded operator on $L^{p'}(\partial\Omega)$, with norm depending on $\Lambda$, $p$, the constants in (\ref{square2}) and the $(D)_{p'}^{A^t}$ constant. Indeed,
\[
\begin{aligned}
& \|H_{A^t}(f_0)\|_{L^{p'}(\partial\Omega)} = 
\|\widetilde{u}\|_{L^{p'}(\partial\Omega)} \leq
\|N(\widetilde{u})\|_{L^{p'}(\partial\Omega)} \\
& \simeq \|S(\widetilde{u})\|_{L^{p'}(\partial\Omega)} \simeq
\|S(u)\|_{L^{p'}(\partial\Omega)}
\simeq \|N(u)\|_{L^{p'}(\partial\Omega)} \leq C(p)\|f_0\|_{L^{p'}(\partial\Omega)},
\end{aligned}
\]
Now let us prove the identity
\[
\int_{\partial\Omega} (\partial_\t(H_Ag)) h - g(\partial_\t(H_{A^t}h)) \,
d\sigma = 0
\]
for any $g, h \in C_0^\infty(\partial\Omega)$. This
is easily done as follows: Let $w$
and $v$ be the extentions of $h$ and $g$ respectively satisfying $\div(A^t
\nabla w) = \div(A \nabla v) = 0$ in $\Omega$ so that $|w(\x)| + |\x||\nabla w(\x)| + |v(\x)| + |\x||\nabla v(\x)| = O(|\x|^{\delta-1})$ for large $|\x|$. Then
\begin{eqnarray*}
& & \int_{\partial\Omega} (\partial_\t(H_Ag)) h - g(\partial_\t(H_{A^t}h)) \,
d\sigma \\
& = & \int_{\partial\Omega} (\n \cdot A \nabla v) h - g(\n \cdot A^t \nabla w) \,
d\sigma \\
& = & \int_{\Omega} \div(A(w\nabla v)) - \div(A^t(v\nabla w)) \, dxdy \\
& = & \int_{\Omega} \nabla w \cdot A\nabla v - \nabla v \cdot A^t\nabla w \, dxdy
= 0.
\end{eqnarray*}
The integrating by parts above in unbounded domains can be justified by the decay properties of $v$ and $w$, passing to the limit through bounded domains.

Combining these two facts we may now prove that
\begin{equation} \label{comp}
\|\n \cdot A\nabla u\|_{L^p(\partial\Omega)} \leq C(p) \|\partial_\t
f_0\|_{L^p(\partial\Omega)},
\end{equation}
for $u$ which solve (\ref{dirichlet}) with boundary data $f_0$. We can do this
using duality:
\begin{eqnarray*}
\left| \int_{\partial\Omega} (\n \cdot A \nabla u)h \, d\sigma \right| & = &
\left| \int_{\partial\Omega}
(\partial_\t(H_{A}f_0))h \, d\sigma \right| \\
& = & \left| -\int_{\partial\Omega}
(\partial_\t f_0) (H_{A^t}h) \, d\sigma \right| \\
& \leq & C(p) \|\partial_\t f_0
\|_{L^p(\partial\Omega)}
\|h\|_{L^{p'}(\partial\Omega)}
\end{eqnarray*}
so (\ref{comp}) follows.
It is a simple
exercise to show that $|\n \cdot \nabla u| \leq C (|\n \cdot A \nabla u| +
|\partial_\t u|)$ from which it is easy to deduce (\ref{comp2}).
We also have the reverse of (\ref{comp}):
\[
\|\partial_\t {u}\|_{L^p(\partial\Omega)} \leq C(p)\|\n \cdot \widetilde{A}\nabla
{u}\|_{L^p(\partial\Omega)}
\]
for $u$ which satisfy (\ref{neumann}) with boundary data $g_0$ and $A$ replaced by $\widetilde{A}$. This
follows
from (\ref{comp}), since the conormal derivative becomes the tangential
derivative of the conjugate
and the tangential derivative becomes the conormal derivative of the
conjugate. As before we then deduce (\ref{comp3}).
\end{proof5}

The following lemma \cite{KKPT} regarding a certain change of variables will be crucial in
obtaining the boundedness of the layer potentials.
\begin{lemma} \label{5}
Suppose $\Omega = \{(x,t) \in \R^2
\, | \, \phi(x) < t \}$ is the domain above the graph of a
Lipschitz function $\phi$. Let $A = A(x)$ be any matrix satisfying the ellipticity
condition (\ref{ellipticity}) and with coefficients independent of the vertical
direction. Also suppose that $\div \, A\nabla u = 0$ in $\Omega$. Then there
exists a change of variables $\Phi \colon \Omega' \to \Omega$ such that
\begin{enumerate}

\item[(i)] If $v = u \circ \Phi$ then $\div \, B \nabla v = 0$ in $\Omega'$,
where $B$ is upper triangluar and independent of the $t$-variable of the
form
\begin{equation} \label{form}
B = \left(  \begin{array}{cc} 1 & c \\ 0 & d \end{array} \right).
\end{equation}

\item[(ii)] The domain $\Omega'$ is the domain above the graph of a Lipschitz
function.

\end{enumerate}

\end{lemma}

\begin{proof}
The transformation $\Phi$ defined by $\Phi(y,s) = (f(y),s + g(y))$ does the job
with $f^{-1}$ chosen to be the primitive of $1/a_{11}$ and $g$ chosen to be
the primitive of $a_{21} \circ f$, where $A = (a_{ij})_{ij}$. We refer the
reader to Lemma 3.47 in \cite{KKPT} for a detailed proof.
\end{proof}

\begin{remark} \label{5.5}
The same result holds with $B$ replaced by the lower triangular matrix $B^t$.
The proof only needs to be modified by choosing $g$ to instead be the primitive
of $a_{12} \circ f$.
\end{remark}

Escauriaza observed (see \cite[p250]{KKPT}) that the operator $\div \,
B^t\nabla$, where $B$ is as in (\ref{form}), can be written in non-divergence form. Indeed, since the coefficients
only depend on the
$x$-direction,
\begin{eqnarray*}
\div \, B^t \nabla u & = & \partial_{xx} u(x,t) + \partial_t(c(x)\partial_xu(x,t)) +
\partial_t(d(x)\partial_tu(x,t)) \\  & = & \partial_{xx} u(x,t) +
c(x)\partial_{xt}u(x,t) +
d(x)\partial_{tt}u(x,t).
\end{eqnarray*}
We remark that then $(u_t, u_x)$ is a conjugate pair, since $u_t$ is a solution
and
we have
\[ 
\begin{array}{rcl}
\partial_t(u_x) & = & \partial_x(u_t) \\
-\partial_x(u_x) & = & c \partial_x(u_t) + d \partial_t(u_t),
\end{array}
\]
so in this case $\widetilde{u_t} = u_x$.
Solutions to non-divergence form equations in $\R^2$ also enjoy additional
regularity. This is stated explicitly as the lemma
below, which is an immediate consequence of Theorem 11.3 in \cite{GT}, Remark
(3) which follows it, and the discussion in Section 11.2. (This can also be seen directly when the coefficients are $t$-independent, since both $u_t$ and $u_x$ solve divergence form equations. Clearly $u_t$ satisfies the same equation as $u$ and $u_x$ solves an equation of the same type since $\widetilde{u_t} = u_x$.)

\begin{lemma} \label{gilbargquasi}
Let $u \in H^{1,2}(B_r)$ satisfy $\div \, B^t \nabla u = 0$ with
coefficient matrix $B^t$, with $B$ as in
(\ref{form}). Then there exists constants $C
> 0$ and $\alpha \in (0,1]$, depending only on $\lambda$ and $\Lambda$, such
that
\[
|\nabla u(\x) - \nabla u(\y)| \leq \frac{C}{r^\alpha} \|\partial_t
u\|_{L^\infty(B_r)} |\x - \y|^\alpha
\]
for all $\x, \y \in B_{2r/3}$.
\end{lemma}

\section{Reduction to the boundedness of layer potentials}

\label{reduct}

The aim of this section is to reduce the proof of Theorem \ref{main} to proving
the boundedness on $L^p(\partial\Omega)$ ($1 < p < \infty$) of the double layer
potential
\[
\mathcal{K}(f)(\x) = \lim_{h \searrow 0} \, \int_{\partial\Omega} \n(\y) \cdot
A^t(\y)\nabla\Gamma_{(x,\phi(x)+h)}(\y)f(\y) \, d\sigma(\y)
\]
and the related potential
\[
\mathcal{L}(f)(\x) = \lim_{h \searrow 0} \, \int_{\partial\Omega} \t(\y) \cdot
\nabla\Gamma_{(x,\phi(x)+h)}(\y)
f(\y) \, d\sigma(\y),
\]
where 
$\x = (x,\phi(x)) \in \partial\Omega$. This will be done in three steps under the assumption that the
coefficient matrix $A$ is of the form (\ref{form}). Firstly we will show
that the non-tangential maximal function of the gradient of a solution can be
controlled in $L^p$-norm by the $L^p$-norms of the layer potentials of $u_t$ and its
conjugate $\widetilde{u_t}$. Secondly, we will show that
$\|\widetilde{u_t}\|_{L^p(\partial\Omega)}
\leq C(p)\|\nabla u\|_{L^p(\partial\Omega)}$. Finally, it is straightforward to combine these results with Theorem \ref{edge} and Lemma \ref{5} to achieve the our aim. 

Lemma \ref{5}
shows us that for the proof of Theorem \ref{main} we may assume $A$ is of the form (\ref{form}) without loss of
generality,
so let us fix, once and for all, $\Gamma_\x$ to be the fundamental solution of
the operator $L = \div \, A \nabla \cdot$ with pole at $\x$, where $A$ is of
this form. With this choice $\Gamma_\x$ is a solution to a
non-divergence form elliptic equation $L^t\Gamma_\x = 0$ away from $\x$. Throughout the next two sections we will make the following a priori assumptions: $A = I$
for large $x$, $A$ and $\phi$ are smooth functions, $\phi' \equiv
\alpha_0$ for $x$ large, and $x \mapsto \phi(x) - \alpha_0x \in C^\infty_0$. Once the theorems here have been proved under our a
priori assumptions it is a simple matter to obtain the general case.

\begin{theorem} \label{inside}
Let $\Omega = \{ (x,t) \in \R^2 \, | \, \phi(x) < t \}$ for
some Lipschitz function $\phi$, $\|\phi'\|_{L^\infty(\R)} \leq k$. Let $L = \div \, A
\nabla$ be an elliptic operator satisfying (\ref{ellipticity}) with coefficient
matrix $A = A(x)$ of measurable functions bounded by $\Lambda$ independent
of the $t$-variable and of the form (\ref{form}). Then for each $p>1$ there
exists a constant $C(p)$, depending only on $\lambda$,
$\Lambda$, $k$ and $p$, such that any function $u \colon \Omega \to \R$ such that $Lu
=0$, and $u(\x) = O(|\x|^{\delta-1})$ and $|\nabla u(\x)| = O(|\x|^{\delta-2})$ for all $\delta > 0$ as $|\x| \to \infty$, we have
\[
\|\widetilde{N}(\nabla{u}) \|_{L^p(\partial\Omega)} \leq C(p)(\|\nabla u\|_{L^p(\partial\Omega)} + \|\mathcal{K}(u_t)\|_{L^p(\partial\Omega)} +
\|\mathcal{L}(\widetilde{u_t})\|_{L^p(\partial\Omega)}).
\]
\end{theorem}

\begin{proof} Firstly, exactly as in the proof of Lemma 8.10 in \cite{KP}, we
will estimate $\widetilde{N}(\nabla u)$ pointwise by the maximal function of $N(u_t)$ and
$\partial_\t u$, so to prove the lemma we will then need to control $N(u_t)$ in $L^p$-norm. This
can be done by using an idea of Verchota and Vogel \cite{VV} to write $u_t$ as
the sum of two potentials.

Fix $\q \in \partial\Omega$ and an $\x \in \Gamma(\q)$. For $\z = (z,r)
\in B_{\frac{3}{4}\delta(\x)}(\x)$, let $\z^* = (z, \phi(z))$ be the vertical projection onto the
boundary of $\Omega$. Choose $m$ so that $(B_{\frac{3}{4}\delta(\x)}(\x))^* \subset
\Delta_{m\delta(\x)}(\q) \equiv \Delta$. For any $c_B \in \R$ Lemmata \ref{cacciopoli}
and \ref{sup} give
\[
\left(\fint_{B(\x,\delta(\x)/2)} |\nabla u|^2\right)^\frac{1}{2} \leq \frac{C}{\delta(\x)}
\fint_{B(\x,3\delta(\x)/4)} |u - c_B|.
\]
Choosing $c_B = \fint_\Delta u$ we have
\[
|u(\z) - \fint_\Delta u| \leq |u(\z) - u(\z^*)| + |u(\z^*) - \fint_\Delta u|.
\]
Now,
\[
\fint_{B(\x,3\delta(\x)/4)} |u(\z^*) - \fint_\Delta u| \leq C \fint_{\Delta}
|u(\z^*) - \fint_\Delta u| \, d\sigma \leq C \delta(\x) \fint_\Delta
|\partial_\t u|
\]
by Poincar\'e's inequality on $\partial\Omega$, where $\partial_\t$ is the
tangential derivative along $\partial\Omega$, and also $|u(\z) - u(\z^*)| \leq
C\delta(\x)N(u_t)(\z^*)$. Therefore, combining these estimates, we have
\[
\widetilde{N}(\nabla u)(\q) \leq C( M(\partial_\t u)(\q) + M(N(u_t))(\q)),
\]
where $M$ is the Hardy-Littlewood maximal operator, and so, for each $p > 1$,
\begin{equation} \label{step1}
\|\widetilde{N}(\nabla u) \|_{L^p(\partial\Omega)} \leq C(p)(\|\nabla u\|_{L^p(\partial\Omega)} + \|N(u_t)\|_{L^p(\partial\Omega)}).
\end{equation}

Now to estimate the last term above, we wish to find a representation for $u_t$. Recall Green's second identity: Let us write $L = \div \, A \nabla$ and
$L^t = \div \, A^t \nabla$, then we have
\[
\iint_\Omega (Lu)v - u(L^tv) = \int_{\partial\Omega} (\n \cdot A\nabla u)v - (\n
\cdot A^t\nabla v)u \, d\sigma
\]
so, for $u$ such that $Lu = 0$ and replacing $v$ with the fundamental solution
$\Gamma_\x$ for $L$, so that
$L^t\Gamma_\x = \delta_\x$, the Dirac mass at $\x$, we obtain
\[
u(\x) = \int_{\partial\Omega} (\n \cdot A^t\nabla \Gamma_\x)u - (\n
\cdot A\nabla u)\Gamma_\x \, d\sigma.
\]
Since the elements of $A$ only
depend on the $x$-variable $u_t$ is also a solution and so we have
\[
u_t(\x) = \int_{\partial\Omega} (\n \cdot A^t\nabla \Gamma_\x)u_t - (\n
\cdot A\nabla u_t)\Gamma_\x \, d\sigma.
\]
Recall that
$\n \cdot A \nabla u_t = \t \cdot \nabla \widetilde{u_t}$ where $\widetilde{u_t}$
is the conjugate of $u_t$. Therefore,
\begin{equation} \label{pot}
\begin{aligned}
u_t(\x) & =  \int_{\partial\Omega} (\n \cdot A^t\nabla \Gamma_\x)u_t - (\partial_\t
(\widetilde{u_t}))\Gamma_\x \, d\sigma \\
& =  \int_{\partial\Omega} (\n \cdot A^t\nabla \Gamma_\x)u_t + 
\t \cdot \nabla \Gamma_\x (\widetilde{u_t}) \, d\sigma. 
\end{aligned}
\end{equation}

To complete the proof we will estimate $N(u_t)$  by estimating each term on
the right-hand side of
(\ref{pot}). This can be done using the singular integral
representation we obtain for $\mathcal{K}$ and $\mathcal{L}$ in Section
\ref{boundedness} and so we leave the end of the proof to be completed in Remark
\ref{endofproof}. It is in this step where we will use the hypothesis that $A$
is of the form (\ref{form}).
\end{proof}

\begin{lemma} \label{transcomp}
Let $u$, $\Omega$ and $A$ be as in Theorem \ref{inside}. There exists a constant $C(p) >
0$ depending only on $\lambda$, $\Lambda$, $k$ and $p$ such that
\[
\|\widetilde{u_t}\|_{L^p(\partial\Omega)} \leq C(p)\|\nabla u\|_{L^p(\partial\Omega)}.
\]
\end{lemma}

\begin{proof}
This would be immediate if it were that $A = B^t$ where $B$ is as in
(\ref{form}) since, for a solution 
$v$ to 
$\div \, B^t\nabla v = 0$,
$(v_t, v_x)$ is a conjugate pair (see Section \ref{prelim}). For our case we can use the change of variables $\Phi$
from Lemma
\ref{5} (and Remark \ref{5.5}), so 
$v = u \circ \Phi$ and 
$B^t = (\det \Phi')(\Phi'^{-1})^t A \Phi'^{-1}$. A simple calculation reveals that taking derivatives in the $t$-variable commutes with the transformation
$\Phi^{-1}$. Also the explicit form of the transformation shows us that
$t$-derivatives of $\widetilde{u}$ are comparable to $t$-derivatives of
$\widetilde{v}$ and $\|\nabla
u\|_{L^p(\partial\Omega)}
\simeq \|\nabla v\|_{L^p(\partial\Omega')}$. Finally we claim conjugation also commutes with the transformation $\Phi^{-1}$. Proving this means showing 
\begin{equation} \label{commute}
\widetilde{v} \circ \Phi^{-1} = (v \circ \Phi^{-1})\widetilde{\hspace{.2cm}}. 
\end{equation}
With $v_0 := \widetilde{u} \circ \Phi$, we have that 
\[
\nabla v = \Phi' (\nabla u \circ \Phi) \hspace{.2cm} \mbox{and} \hspace{.2cm} \nabla v_0 = \Phi' (\nabla \widetilde{u} \circ \Phi),
\]
so using (\ref{conj}), we obtain
\[
\left(  \begin{array}{cc} 0 & 1 \\ -1 & 0 \end{array} \right)(\Phi')^{-1} \nabla v_0 = A(\Phi')^{-1} \nabla v.
\]
Multiplying on the left by 
$(\Phi'^{-1})^t$ we see that
\[
(\det\Phi')^{-1}\left(  \begin{array}{cc} 0 & 1 \\ -1 & 0 \end{array} \right) \nabla v_0 = (\Phi'^{-1})^t A (\Phi')^{-1} \nabla v.
\]
Thus $(v,v_0)$ solves (\ref{conj}) with 
$A$ replaced with $B^t$, that is
$v_0 = \widetilde{v}$. The last equality can be rewritten 
$\widetilde{u} \circ \Phi = (u \circ \Phi)\widetilde{\hspace{.2cm}}$ from which (\ref{commute}) follows.

These observations allow us to conclude that
\[
\begin{aligned}
&\| \widetilde{u_t}\|_{L^p(\partial\Omega)} = \| (\partial_t(v \circ \Phi^{-1}))\widetilde{\hspace{.2cm}}\|_{L^p(\partial\Omega)}
=
\|\partial_t(\widetilde{v} \circ \Phi^{-1})\|_{L^p(\partial\Omega)} \\
& \leq C(p)\|\partial_t \widetilde{v}\|_{L^p(\partial\Omega')} =
C(p)\|\widetilde{v_t} \|_{L^p(\partial\Omega')} \leq C(p)\|\nabla
v\|_{L^p(\partial\Omega')} \simeq \|\nabla u\|_{L^p(\partial\Omega)},
\end{aligned}
\]
and so the lemma is proved.
\end{proof}

Once we know $\mathcal{K}$ and $\mathcal{L}$ are bounded for all $A$ of the form (\ref{form}), we have from Theorems \ref{inside} and \ref{transcomp} that $\|\widetilde{N}(\nabla u)\|_{L^p(\partial\Omega)} \leq C(p)\|\nabla u\|_{L^p(\partial\Omega)}$ ($1<p<\infty$). Using Lemma \ref{5} we extend this to arbitrary $A$. One then applies Theorems \ref{KKPT} and \ref{edge} to conclude the proof of Theorem \ref{main}. It remains to demonstrate the boundedness of $\mathcal{K}$ and $\mathcal{L}$ on $L^p$. This task will occupy us for the next two sections.

\section{Boundedness of the layer potentials on boundaries with small Lipschitz
constants}

\label{boundedness}

In this section we will show that $\mathcal{K}$ and $\mathcal{L}$ are bounded on
$L^p(\partial\Omega)$ ($1 < p < \infty$) in the special case that $\Gamma_\x$ is the
fundamental solution of the operator $L = \div(A \nabla \, \cdot)$ with pole at
$\x$, where $A$ is of the form (\ref{form}) and that our domain $\Omega$ has a
boundary which is close to linear in the sense that, for $\alpha_0 \in
[-k,k]$, $\|\phi' - \alpha_0
\|_{L^\infty(\R)} \leq \varepsilon_0$, where $\varepsilon_0 > 0$ depends only on the
ellipticity constants and $k>0$. Lemma \ref{5} shows that the first assumption can be made
without loss of generality. However, if we are to succeed in proving Theorem
\ref{main}, the second assumption must be dropped. This will be done in Section
\ref{david}. Denote by $\Lambda^k(\varepsilon_0)$ the set of all Lipschitz
functions $\phi$ such that $\|\phi' - \alpha_0
\|_{L^\infty(\R)} \leq \varepsilon_0$, with $\alpha_0 \in
[-k,k]$. We will always require that $0 < \varepsilon_0 \leq k$ so the Lipschitz
constant of $\phi \in \Lambda^k(\varepsilon_0)$ is no more than $2k$. We remind the reader of the a priori smoothness assumptions made at the beginning of Section~\ref{reduct}.

\begin{lemma} \label{kernelest1}
There
exists constants $C_3 > 0$ and $0 < \alpha \leq 1$, depending only on $\lambda$
and $\Lambda$, such that
\begin{equation} \label{103}
|\nabla\Gamma_\x(\y)| \leq \frac{C_3}{|\x-\y|}
\end{equation}
for all $\x,\y \in \R^2$, $\x \neq \y$,
\begin{equation} \label{104}
|\nabla\Gamma_\x(\y) - \nabla\Gamma_{\x'}(\y)| \leq \frac{C_3|\x-\x'|^\alpha}{(|\x-\y| +
|\x'-\y|)^{1 + \alpha}}
\end{equation} 
when $|\x-\x'| \leq \frac{1}{2}\max\{|\x-\y|,|\x'-\y|\}$, and
\begin{equation} \label{105}
|\nabla\Gamma_\x(\y) - \nabla\Gamma_\x(\y')| \leq \frac{C_3|\y-\y'|^\alpha}{(|\x-\y| + |\x-\y'|)^{1 + \alpha}}
\end{equation}
when $|\y-\y'| \leq \frac{1}{2}\max\{|\x-\y|,|\x-\y'|\}$.
\end{lemma}

\begin{proof}
Fix $\x,\y \in \R^2$, $\x \neq \y$ and set $r = |\x - \y|$. Cacciopoli's
inequality (Lemma \ref{cacciopoli}) gives us that
\begin{equation} \label{101}
\begin{aligned}
\fint_{B_{r/4}(\y)} |\nabla \Gamma_{\x}(\z)|^2 \, d\z & \leq \frac{C}{r^2} \,
\fint_{B_{r/2}(\y)} |\Gamma_\x(\z)|^2 \, d\z \\
\end{aligned}
\end{equation}
Now we apply a linear change of variables $\Phi$
which sends $\0 = (0,0)$ to $\x$ and $(0,-R_1/4)$ to $\y$:
\begin{equation} \label{101.5}
\begin{aligned}
& \frac{C}{r^2} \,
\fint_{B_{\frac{r}{2}}(\y)} |\Gamma_\x(\z)|^2 \, d\z = \frac{C}{r^2} \,
\fint_{B_{\frac{1}{8}R_1}((0,-R_1/4))} |\Gamma_{\Phi(\0)} \circ \Phi (\z)|^2 \, d\z\\
& \leq \frac{C}{r^2} \,
\fint_{B_{\frac{1}{8}R_1}((0,-R_1/4))} |\log(|\z|)|^2 \, d\z \leq \frac{C}{r^2},
\end{aligned}
\end{equation}
where the second to last inequality follows from Theorem \ref{fundamental},
since $\Gamma_{\Phi(\0)} \circ
\Phi$ is the fundamental solution to an elliptic operator with the same
constants $\lambda$ and $\Lambda$ and pole at the origin.
Now $\partial_s \Gamma_\x(y,s)$ satisfies
$L(\partial_s \Gamma_\x)
= 0$ away from $\x$ since the coefficients do not depend on the second variable, thus, using
the above and Lemma
\ref{sup}, we have
\begin{equation} \label{102}
\sup_{B_{r/8}(\y)}|\partial_s \Gamma_\x| \leq
\left(\fint_{B_{r/4}(\y)} |\nabla \Gamma_{\x}(\z)|^2 \, d\z \right)^\frac{1}{2} \leq \frac{C}{r}
\end{equation}
This last bound shows we may now apply Lemma \ref{gilbargquasi} with
$B_r$ replaced with $B_{r/8}(\y)$ to obtain
\begin{eqnarray*}
|\partial_y \Gamma_\x(\y)| & \leq & |\partial_y \Gamma_\x(\z)| + |\partial_y
\Gamma_\x(\y) - \partial_y \Gamma_\x(\z)| \\
& \leq & |\partial_y \Gamma_\x(\z)|
+ \frac{C}{r} \, \left(\frac{|\y - \z|}{r}\right)^\alpha \\
& \leq & |\partial_y \Gamma_\x(\z)|
+ C/r
\end{eqnarray*}
for $\z \in B_{r/8}(\y)$. Averaging this inequality over $\z \in B_{r/8}(\y)$
and using (\ref{101}) and (\ref{101.5}) with $r$ replaced with $r/2$ yields
\[
|\partial_y \Gamma_\x(\y)| \leq \frac{C}{r}.
\]
Combining this with (\ref{102}) yields (\ref{103}). The same application of
Lemma \ref{gilbargquasi} now easily yields (\ref{105}): with $r = \max(|\y -
\x|, |\y' - \x|) = |\y - \x|$, say, we have for $|\y - \y'| \leq r/2$ that
\[
\begin{aligned}
|\nabla \Gamma_\x(\y) - \nabla \Gamma_\x(\y')| & \leq \frac{C}{r^\alpha} \|\nabla
\Gamma_\x \|_{L^\infty(B_{3r/4}(\y))} |\y - \y'|^\alpha \\
& \leq \frac{C}{r^{1 + \alpha}}|\y - \y'|^\alpha.
\end{aligned}
\]
Inequality (\ref{104}) follows
from the fact (Lemma \ref{symmetry}) that $\x \mapsto
\Gamma_\x(\y)$ is a solution to an elliptic equation in divergence form, so we may use the
standard H\"older continuity result of Lemma \ref{Calpha}: now with $r =
\max(|\x -
\y|, |\x' - \y|) = |\x - \y|$, say, we have for $|\x - \x'| \leq r/2$ that
\[
\begin{aligned}
|\nabla \Gamma_\x(\y) - \nabla \Gamma_{\x'}(\y)| & \leq \frac{C}{r^\alpha} \|\nabla
\Gamma_\x \|_{L^2(B_{3r/4}(\y))} |\x - \x'|^\alpha \\
& \leq \frac{C}{r^{1 + \alpha}}|\x - \x'|^\alpha,
\end{aligned}
\]
which completes the proof of the lemma.
\end{proof}

Let $\phi \in \Lambda^\frac{k}{4}(\varepsilon_0)$ be a Lipschitz function with $\|\phi' -
\alpha_0\|_{L^\infty(\R)} = \varepsilon \leq \varepsilon_0$, $\Omega = \Omega^+ = \{
(x,t) \, | \, t > \phi(x) \}$ and $\Omega^- = \{
(x,t) \, | \, t < \phi(x) \}$, and write $\partial\Omega_h = \partial\Omega + (0,h)$
and
\[
\left(  \begin{array}{c} k_{1,h}(x,y) \\ k_{2,h}(x,y) \end{array} \right) =
\nabla\Gamma_{(x,\phi(x)+h)}(y,\phi(y))
\]
for $h \in \R$. Define the function $K_h \colon \R^2 \to
\mathcal{M}$ by
\begin{equation} \label{kernel}
K_h = \left(  \begin{array}{cc} k_{1,h} & k_{2,h} \\ k_{1,h} & k_{2,h}
\end{array} \right).
\end{equation}
We will write $K_0 = K$. We have the following easy corollary to Lemma
\ref{kernelest1}
\begin{corollary} \label{kernelest2}
For each $\phi \in \Lambda^\frac{k}{4}(\varepsilon_0)$, the function $K \colon \R^2 \to \M$ is a Calder\'on-Zygmund kernel. That is to say, there
exists constants $C_4 = C_4(k) > 0$ and $\alpha \in (0,1]$, depending only on $\lambda$,
$\Lambda$ and the Lipschitz constant $k$, such that
(\ref{kerest1}), (\ref{kerest2}) and (\ref{kerest3}) hold with $C=C_4$.
\end{corollary}

For Banach spaces $\mathcal{X}$ and $\mathcal{Y}$ of real matrix-valued
functions on $\R$,
define $T \colon \mathcal{X} \to \mathcal{Y}'$ by
\begin{equation} \label{def}
\langle G,T(F)\rangle = \lim_{h \searrow 0} \iint_{\R^2} G(x)^tK_h(x,y)F(y)
\, dydx,
\end{equation}
where $F \in \mathcal{X}$ and $G \in \mathcal{Y}$ and $\langle \cdot, \cdot\rangle$ is the dual
pairing in $\mathcal{Y}$.
The transpose operator $T^t$ is defined as $\langle F,
T^t(G)\rangle = \langle G, T(F)\rangle^t$. Our aim will be to show that $T$
extends to a
bounded operator on
$L^2(\R^2,\mathcal{M})$, but first we must choose appropriate Banach spaces,
$\mathcal{X}$ and $\mathcal{Y}$, and show this definition makes sense. This is the content of the
following lemma, and for that we will need the following notation.
We denote by $\S = \S(\R,\M)$ the
space of Schwarz functions on $\R$ with values in $\M$ and for a function $B_0
\colon \R \to \M$, $B_0\S$ is the set of all functions obtained by
multiplying functions in $\S$ on the left by $B_0$. Let $\delta$ be the
distance
function $\delta(\x) = \mbox{dist}(\x,\partial\Omega)$. The spaces
$C^\infty_0(\R,\R)$ and $C^\infty_0 = C^\infty_0(\R,\M)$ are the spaces of compactly
supported smooth functions from $\R$ to $\R$ or $\M$ respectively.
Set $B_1 \colon \R \to \M$ equal to the matrix with first column
being $(1+(\phi')^2)^\frac{1}{2}A\n$ and second column $(1+(\phi')^2)^\frac{1}{2}\t$.
Set $B_2 \colon \R \to \M$ equal to the diagonal matrix with both
diagonal entries being $(1+(\phi')^2)^\frac{1}{2}\n \cdot A^t \kappa^\perp$, where
$\kappa^\perp = (-\alpha_0, 1)$, and $B_3 \colon \R \to \M$ equal to the diagonal matrix with both diagonal
entries being $(1+(\phi')^2)^\frac{1}{2}\t \cdot \kappa$, where $\kappa =
(1,\alpha_0)$.

Note that by the ellipticity condition (\ref{ellipticity}), for sufficiently small $\varepsilon_0$
(depending only on the ellipticity constants and $k$), $B_1$, $B_2$ and $B_3$
are invertible with inverses bounded in $L^\infty(\R, \M)$-norm in terms of the
ellipticity constants and $k$.

\begin{lemma} \label{continuous}
For each 
$\phi \in \Lambda^\frac{k}{4}(\varepsilon_0)$, the operator $T$ is a continuous linear operator from $B_1\S$ to $(B_0\S)'$, for
any bounded $B_0$.
\end{lemma}

\begin{proof}
We fix $f, g \in \S(\R, \R)$, $h \in C^\infty_0(\R,\R)$ with $h$ positive and
equal to one near zero and define $u,v \colon \Omega \to \R$ by
\[
u(\x) = f(x)h(t - \phi(x)) \hspace{0.2cm} \mbox{and} \hspace{0.2cm} v(\x) =
g(x)h(t - \phi(x)).
\]
Green's first identity gives us that
\[
u(\x) = -\iint_\Omega \nabla u(\y) \cdot A^t(\y)\nabla \Gamma_\x(\y) \, d\y + 
\int_{\partial\Omega} u(\y) \n(\y) \cdot A^t(\y)\nabla \Gamma_\x(\y) \, d\sigma(\y).
\]
By multiplying by a bounded function $b$ and the function $v$, and integrating,
we obtain
\begin{equation} \label{cont2}
\begin{aligned}
& \int_{\partial\Omega_h} u(\x)b(\x)v(\x) \, d\x =\\
& \int_{\partial\Omega_h} \left(-\iint_\Omega \nabla u(\y) \cdot A^t(\y)\nabla \Gamma_\x(\y) \, 
d\y\right) b(\x)v(\x) \, d\sigma(\x)  \\
&+ \int_{\partial\Omega_h} \int_{\partial\Omega} u(\y) \n(\y) \cdot A^t(\y)\nabla
\Gamma_\x(\y) b(\x)v(\x)\, d\sigma(\y) \, d\sigma(\x).
\end{aligned}
\end{equation}
The left-hand side can easily be controlled by the bound on $b$ and the product
of semi-norms $[f]^{0}_{1}[g]^{0}_{1}$, where
\begin{equation} \label{norms}
[f]^\beta_N = \sup_\x (1+ |\x|)^N|\partial^\beta f(\x)|
\end{equation}
for $N \in {\bf{Z}}$ and $\beta \in \N$. We can control the first term on the right-hand side similarly.
Indeed, observe that
\[
\begin{aligned}
\left| \iint_\Omega \nabla u(\y) \cdot A^t(\y)\nabla \Gamma_\x(\y) \, 
d\y \right| & \leq C([f]^{0}_3 + [f]^1_3) \iint_{\R^2} \frac{(1+|\y|)^{-3}}{|\x
- \y|} \, d\y \\
& \leq C([f]^{0}_3 + [f]^1_3)
\end{aligned}
\]
and then proceed as before. And so, we conclude the remaining term on the
right-hand side is also controlled by a finite sum of products
$\|b\|_{L^\infty(\R)}[f]^\beta_N[g]^\gamma_M$, that is
\begin{equation} \label{cont}
\begin{aligned}
\left| \int_{\partial\Omega_h} \int_{\partial\Omega} u(\y) \n(\y) \cdot  A^t(\y) \nabla
\Gamma_\x(\y) b(\x)v(\x)\, d\sigma(\y) \, d\sigma(\x) \right| \\
\leq C\|b\|_{L^\infty(\R)}\sum_{N,M,
\beta, \gamma}
[f]^\beta_N[g]^\gamma_M.
\end{aligned}
\end{equation}
We may also compute, setting $R_3 = R_1/(1 + k^2)^\frac{1}{2}$
(where $R_1$ is the constant appearing in Theorem \ref{fundamental}),
\[
\begin{aligned}
& \left| \int_{\partial\Omega} (\t(\y) \cdot \nabla \Gamma_\x(\y)) u(\y) \,
d\sigma(\y) \right| \\
= & \left| \int_{\R} \{\partial_1\Gamma_{(x,\phi(x))}(y,\phi(y)) +
\phi'(y)\partial_2\Gamma_{(x,\phi(x))}(y,\phi(y))\} u(\y) \, dy \right| \\
= & \left| \int_{\R} \frac{d}{dy}\{\Gamma_{(x,\phi(x))}(y,\phi(y))\} u(y,\phi(y))
\, dy \right| \\
\leq & \left| \int_{\R \setminus (x - R_3,x + R_3)}
\frac{d}{dy}\{\Gamma_{(x,\phi(x))}(y,\phi(y))\} u((y,\phi(y))
\, dy \right| \\
& + \left| \int_{(x - R_3,x + R_3)} \Gamma_{(x,\phi(x))}(y,\phi(y)) \frac{d}{dy}u(y,\phi(y))
\, dy \right| \\
& + \sum_{\y \in \{x - R_3,x + R_3\}}
|\Gamma_{(x,\phi(x))}(y,\phi(y)) u(y,\phi(y))| \\
\leq & C([f]_2^{0} + [f]_0^{1} + [f]_0^{0}),
\end{aligned}
\]
so allowing us to conclude
\begin{equation} \label{cont3}
\begin{aligned}
\left| \int_{\partial\Omega_h} \int_{\partial\Omega} u(\y) \t(\y) \cdot \nabla
\Gamma_\x(\y) b(\x)v(\x)\, d\sigma(\y) \, d\sigma(\x) \right| \\
 \leq C\|b\|_\infty\sum_{N,M,
\beta, \gamma}
[f]^\beta_N[g]^\gamma_M.
\end{aligned}
\end{equation}

However,
\[
K_h(x,\cdot)B_1 =
\left(  \begin{array}{cc}
(1+\phi')^\frac{1}{2} \n \cdot A^t \nabla\Gamma_{(x,\phi(x)+h)} & 
(1+\phi')^\frac{1}{2} \t \cdot \nabla\Gamma_{(x,\phi(x)+h)}  \\
(1+\phi')^\frac{1}{2} \n \cdot A^t \nabla\Gamma_{(x,\phi(x)+h)} & 
(1+\phi')^\frac{1}{2}\t \cdot \nabla\Gamma_{(x,\phi(x)+h)}
\end{array} \right),
\]
so by writing in (\ref{def})
\[
G(x)^tK_h(x,y)F(y) = G_\infty(x)^t B_0(x)^t K_h(x,y) B_1(y)
F_\infty(y)
\]
with $F_\infty \in \S$ and $G_\infty \in \S$, we see that
each matrix element in (\ref{def}) is a sum of terms of the form of either
(\ref{cont}) or (\ref{cont3}), thus provided the limit exists, (\ref{def}) is controlled by a finite sum of $[F_\infty]^\beta_N[G_\infty]^\gamma_M$, with the semi-norms defined analogously to
(\ref{norms}). We can easily show the limit in (\ref{def}) exists almost everywhere under our a priori assumptions: By the above it suffices to show both
\begin{equation} \label{limit1}
\int_{\partial\Omega} \n(\y) \cdot A^t(\y) \nabla \Gamma_\x(\y) u(\y) \, d\sigma(\y)
\end{equation}
and
\begin{equation} \label{limit2}
\int_{\partial\Omega} \t(\y) \cdot \nabla \Gamma_\x(\y) u(\y) \, d\sigma(\y)
\end{equation}
converge uniformly in 
$x$ as 
$h \searrow 0$, where 
$\x = (x,\phi(x) + h)$. To show (\ref{limit1}) converges we write
\[
\begin{aligned}
& \int_{\partial\Omega} \n(\y) \cdot A^t(\y) \nabla \Gamma_\x(\y) u(\y) \, d\sigma(\y) \\
& = \int_{\partial\Omega \cap B_1(\x)} \n(\y) \cdot A^t(\y) \nabla \Gamma_\x(\y) u(\y) \, d\sigma(\y) \\
& + \int_{\partial\Omega \setminus B_1(\x)} \n(\y) \cdot A^t(\y) \nabla \Gamma_\x(\y) u(\y) \, d\sigma(\y).
\end{aligned}
\]
The second integral on the right-hand side clearly converges uniformly in 
$x$. To show that the first also does, we can use the divergence theorem to see it is equal to
\[
\begin{aligned}
& \int_{\partial\Omega \cap B_1(\x)} \n(\y) \cdot A^t(\y) \nabla \Gamma_\x(\y) (u(\y) - u(\x)) \, d\sigma(\y) \\
& + u(\x) \int_{\partial\Omega \cap B_1(\x)} \n(\y) \cdot A^t(\y) \nabla \Gamma_\x(\y) \, d\sigma(\y) \\
& = \int_{\partial\Omega \cap B_1(\x)} \n(\y) \cdot A^t(\y) \nabla \Gamma_\x(\y) (u(\y) - u(\x)) \, d\sigma(\y) \\
& - u(\x)\int_{\Omega^- \cap \partial B_1(\x)} \n(\y) \cdot A^t(\y) \nabla \Gamma_\x(\y) \, d\sigma(\y),
\end{aligned}
\]
and both integrals on the right-hand side converge uniformly in 
$x$ under our a priori smoothness assumptions. We can see (\ref{limit2}) converges similarly: again we write
\[
\begin{aligned}
\int_{\partial\Omega} \t(\y) \cdot \nabla \Gamma_\x(\y) u(\y) \, d\sigma(\y)
& = \int_{\partial\Omega \cap B_1(\x)} \t(\y) \cdot \nabla \Gamma_\x(\y) u(\y) \, d\sigma(\y) \\
& + \int_{\partial\Omega \setminus B_1(\x)} \t(\y) \cdot \nabla \Gamma_\x(\y) u(\y) \, d\sigma(\y).
\end{aligned}
\]
As before the second integral on the right-hand side converges uniformly in 
$x$. The Fundamental Theorem of Calculus shows that the first is equal to
\[
\begin{aligned}
& \int_{\partial\Omega \cap B_1(\x)} \t(\y) \cdot \nabla \Gamma_\x(\y) (u(\y) - u(\x)) \, d\sigma(\y) \\
& + u(\x) \int_{\partial\Omega \cap B_1(\x)} \t(\y) \cdot \nabla \Gamma_\x(\y) \, d\sigma(\y) \\
& = \int_{\partial\Omega \cap B_1(\x)} \t(\y) \cdot \nabla \Gamma_\x(\y) (u(\y) - u(\x)) \, d\sigma(\y) \\
& - u(\x) \int_{\Omega^- \cap \partial B_1(\x)} \t(\y) \cdot \nabla \Gamma_\x(\y) \, d\sigma(\y),
\end{aligned}
\]
and once again both integrals on the right-hand side converge uniformly in 
$x$.
\end{proof}

Since we clearly have the representation (\ref{singrep}), what we have shown so
far is that $T$ is a singular integral operator associated to the
Calder\'on-Zygmund kernel $K$.

\begin{remark} \label{endofproof}
This singular integral representation can be used to complete the proof of Lemma
\ref{inside} since $\mathcal{K}$ and $\mathcal{L}$ can be easily written in
terms of $T$ and so using (\ref{pot}) can control $N(u_t)$.
We define
\[
\mathcal{T}(F)(\z) = \int_{\R} K_r(z,y)F(y) \, dy,
\]
for $(z,r) = \z \in \Omega$, and
\[
T^{(\delta)}(F)(\x) = \int_{|y -x| \geq \delta} K(x,y)F(y) \, dy,
\]
for $\x \in \partial\Omega$. Fixing $(z,r) = \z \in \Gamma(\x)$ we claim that
$|T^{(\delta)}(F)(\x) - \mathcal{T}(\z)| \leq
CM(F)(\x)$ for $\delta = \delta(\z)$. We have
\[
\begin{aligned}
& |T^{(\delta)}(F)(\x) - \mathcal{T}(\z)| \\
& = \left| \int_{|x-y|\geq \delta(\z)} K(x,y)F(y) \, dy -
\int_{\R} K_r(z,y)F(y) \, dy \right| \\
& = \left| \int_{|x-y|\geq \delta(\z)} (K(x,y) - K_r(z,y))F(y) \, dy -
\int_{|x-y| < \delta(\z)} K_r(z,y)F(y) \, dy \right| \\
& \leq \left| \int_{|x-y|\geq \delta(\z)} (K(x,y) - K_r(z,y))F(y) \, dy \right| +
\left| \int_{|x-y| < \delta(\z)} K_r(z,y)F(y) \, dy \right|.
\end{aligned}
\]
The first term can be estimated using (\ref{104}):
\[
\begin{aligned}
&\left| \int_{|x-y|\geq \delta(\z)} (K(x,y) - K_r(z,y))F(y) \, dy \right|
\leq C \int_{|x-y|\geq \delta(\z)} \frac{|x - z|^\alpha |F(y)|}{|x - y|^{1 +
\alpha}} \, dy \\
&\leq C(1 + a)^\alpha\delta(\z)^\alpha \int_{|x-y|\geq \delta(\z)}
\frac{|F(y)|}{|x - y|^{1 +
\alpha}} \, dy
\leq CM(F)(x).
\end{aligned}
\]
And the second term is also easily controlled using (\ref{103}):
\[
\left| \int_{|x-y| < \delta(\z)} K_r(z,y)F(y) \, dy \right|
\leq \frac{C}{\delta(\z)} \int_{|x-y| < \delta(\z)} |F(y)| \, dy \leq
CM(F)(x),
\]
which proves the claim. The claim tells us that
$T^{(\delta(\z))}(F)(\x)$ and $M(F)(\x)$ control
$\mathcal{T}(F)(\z)$, so $N(\mathcal{T})$ can be controlled
in operator norm by the operator norm of $T^* + M$. Standard
Calder\'on-Zygmund theory (see, for example,\cite{G}) tells us that we can control $T^*$ in terms of $T$. Writing $\mathcal{K}$ and $\mathcal{L}$ in terms of $T$ and
using (\ref{pot}) then completes the proof of Lemma \ref{inside}.
\end{remark}

Our aim now is to prove the following.

\begin{theorem} \label{T(B)}
For each $k > 0$, there exists an $\varepsilon_0 > 0$ depending only on $k$,
$\lambda$ and $\Lambda$, such that, for any $\phi \in \Lambda^\frac{k}{4}(\varepsilon_0)$, the singular
integral operator $T$ admits a continuous extension to
$L^2(\R,\M)$ and therefore also to $L^p(\R,\M)$ for all $1<p<\infty$ with norm
depending only on $p$, $\lambda$, $\Lambda$ and $k$. Consequently the
layer potentials $\mathcal{K}$ and $\mathcal{L}$,
defined, for 
$\x = (x,\phi(x)) \in \partial\Omega$, by
\[
\mathcal{K}(f)(\x) = \lim_{h \searrow 0} \, \int_{\partial\Omega} \n \cdot
A^t(\y)\nabla\Gamma_{(x,\phi(x)+h)}(\y)f(\y) \, d\sigma(\y)
\]
and
\[
\mathcal{L}(f)(\x) = \lim_{h \searrow 0} \, \int_{\partial\Omega} \t \cdot
\nabla\Gamma_{(x,\phi(x)+h)}(\y)
f(\y) \, d\sigma(\y),
\]
are bounded linear operators on $L^p(\partial\Omega,\R)$
($1<p<\infty$) when $\Omega = \{\x = (x,t) \, | \, t> \phi(x) \}$.
\end{theorem}



This will be proved by applying a version of the $T(B)$-Theorem. Let us first
state some standard definitions and fix some more notation.

Let $M_{B_0}$ denote the operator which multiplies on the left by a function $B_0 \colon \R
\to \M$. By a \emph{normalised bump} we mean a function $F = (f_{ij})_{ij} \in
\S$ with support
in a ball with radius $10$ such that
\[
[F]^\beta_0 := \sup_{i,j,x} |\partial^\beta f_{ij}(x)| \leq 1
\]
for $\beta = 0,1,2$. For such an $F$ define $F_R = R^{-1}F(\cdot/R)$. An
operator $T$ is said to satisfy the \emph{weak boundedness property} if
there exists a constant $C > 0$ such that
\begin{equation} \label{weakboundednessdef}
|\langle G_R, TF_R \rangle| \leq CR^{-1}.
\end{equation}

A function $f \colon \R \to \R$ is said to be of \emph{bounded mean
oscillation}, written $f \in \mbox{BMO}$, if
we have
\[
\|f\|_{\mbox{\tiny BMO}} := \sup_I \fint_I |f(x) - \fint_I f| \, dx < \infty,
\]
where the supremum is taken over all intervals $I$.
A function $a \colon \R \to \R$ is said to be an \emph{$H^1$ atom} if there
exists an interval $I \subset \R$ such that (i) $a$ is supported in $I$, (ii)
$\|a\|_{L^2(\R)} \leq |I|^{-\frac{1}{2}}$, and (iii) $\int a = 0$.
For a bounded function $B_0 \in B_1C^\infty_0$, $T(B_0)$ may be interpreted as a
matrix-valued function, whose matrix elements are functions in BMO in
the following manner. We shall call a matrix-valued function
whose matrix elements are $H^1$ atoms, or BMO functions, again an $H^1$ atom,
or a BMO
function, respecively.
For each $H^1$ atom $A_0 \in B_2C^\infty_0$ with support in a ball $B$, choose
a smooth function $\eta \colon \R \to \R$ with compact support which is identically one on the
double of the ball
$2B$. Then we define
\begin{equation} \label{defBMO}
\langle A_0, T(B_0) \rangle := \langle A_0, T(\eta B_0) \rangle + \int_\R \left[ \int_\R
A_0(x)K(x,y) \, dx \right] (1-\eta(y))B_0(y) \, dy.
\end{equation}
It is easy to show this definition is independent of the choice of $\eta$ and
each term on the right-hand side is well-defined. Since our a priori assumptions mean that $B_2$ is smooth and has a smooth inverse, the decomposition of 
$H^1$ via smooth atoms and duality determine
$T(B_0)$ as a BMO function.

We remark that, by duality, $T(F)$ can be identified as an element in
$L^\infty(\R,
\M)$ if we can show that $|\langle G, T(F) \rangle| \leq C\|G\|_{L^1(\R,\M)}$, and it
certainly suffices to show that
\[
\sup_{x\in\R} \sup_{0< h < 1} \left| \int_\R K_h(x,y)F(y) \, dy \right| \leq C.
\]
The version of the $T(B)$-Theorem that we will use is the following \cite{DJS}.

\begin{theorem} \label{T(B)David}
Let $B_1$ and $B_2$ be two bounded functions from $\R$ to $\M$, such that their
inverses exist and are bounded, all of these bounds being no more than some
constant $C_5 > 0$.
Let $T \colon B_1\S \to (B_2\S)'$ be a singular integral
operator such that
$M_{B_2^t}TM_{B_1}$ is weakly bounded, with the constants $C$ appearing in (\ref{kerest1}), (\ref{kerest2}), (\ref{kerest3}) and (\ref{weakboundednessdef}) being no more than $C_6 > 0$, and $T(B_1)$ and
$T^t(B_2)$ are in BMO with norms 
no more than $C_7 > 0$. Then $T$ admits a continuous extension onto
$L^2(\R,\M)$ with norm $\|T\|_{L^2 \to L^2} \leq C(C_5 + C_6 + C_7)$.
\end{theorem}

With this result at hand, it is easy to see that Theorem \ref{T(B)} will be a
consequence of the following lemmata.

\begin{lemma}
For each 
$\phi \in \Lambda^\frac{k}{4}(\varepsilon_0)$, the operator $M_{B_0^t}TM_{B_1}$ satisfies the weak boundedness property for any
bounded $B_0$, the constant 
$C$ in (\ref{weakboundednessdef}) depending only on 
$\lambda$, 
$\Lambda$,
$k$ and the bound on 
$B_0$. Also
$T(B_1) \in L^\infty(\R,\M)$, with 
$L^\infty$-norm bounded in terms of $\lambda$, 
$\Lambda$ and
$k$.
\label{weakboundedness}
\end{lemma}

\begin{proof}
We will first prove the weak boundedness of $M_{B_0^t}TM_{B_1}$. For any $i,j = 1,2$ we choose
\begin{equation} \label{tilde2}
\begin{aligned}
u(\x) & =
R^{-1}f_{ij}(\x/R)h((t-\phi(x))/R), \\
v(\x) & = R^{-1}g_{ij}(x/R)h((t-\phi(x))/R)
\end{aligned}
\end{equation}
and
$b(\x) = b_0^{ij}(\phi(x))h(t)$ in (\ref{cont2}), where $F = (f_{ij})_{ij}$ and
$G = (g_{ij})_{ij}$ are normalised bumps, $B_0 =
(b_0^{ij})_{ij}$ is a bounded function, and $h$ is a real-valued normalised bump centred at the origin and
identically equal to one near the origin. We
have
\[
\left|\iint_\Omega \nabla u(\y) \cdot A^t(\y)\nabla \Gamma_\x(\y) \, 
d\y \right| \leq \|\nabla u\|_{L^\infty(\R^2)} \|\nabla \Gamma_\x\|_{L^1(B_{10R})}
\leq CR^{-2+1},
\]
so
\[
\left| \int_{\partial\Omega_h} \left(\iint_\Omega \nabla u(\y) \cdot A^t(\y)\nabla \Gamma_\x(\y) \, 
d\y\right) b(\x)v(\x) \, d\sigma(\x) \right| \leq CR^{-1}.
\]
Also
\[
\left|\int_{\partial\Omega_h} u(\x)b(\x)v(\x) \, d\x \right| \leq
\|b\|_{L^\infty(\R^2)}\|u\|_{L^2(\partial\Omega)} \|v\|_{L^2(\partial\Omega)}
\leq CR^{-1},
\]
so (\ref{cont2}) allows us to conclude
\begin{equation} \label{wbone}
\left| \int_{\partial\Omega_h} \int_{\partial\Omega} u(\y) \n(\y) \cdot A^t(\y)\nabla
\Gamma_\x(\y) b(\x)v(\x)\, d\sigma(\y) \, d\sigma(\x) \right| \leq \frac{C}{R}.
\end{equation}
Now, if the supports of $u$ and $v$ are separated by a distance $R$, then we
have
\begin{equation} \label{wbthree}
\begin{aligned}
& \left| \int_{\partial\Omega_h} \int_{\partial\Omega} (\t(\y) \cdot \nabla
\Gamma_\x(\y)) u(\y) b(\x)v(\x) \,
d\sigma(\y) \, d\sigma(\x) \right| \\
& =
CR^{-1} \, \int_{\partial\Omega_h} \int_{\partial\Omega} |u(\y) b(\x)v(\x)| \,
d\sigma(\y) \, d\sigma(\x) \leq CR^{-1}.
\end{aligned}
\end{equation}
Otherwise, integration by parts shows us that
\[
\begin{aligned}
& \int_{\partial\Omega_h} \int_{\partial\Omega} (\t(\y) \cdot \nabla
\Gamma_\x(\y)) u(\y) b(\x)v(\x) \,
d\sigma(\y) \, d\sigma(\x) \\
& =
- \int_{\partial\Omega_h} \int_{\partial\Omega} \Gamma_\x(\y) (\t(\y) \cdot
\nabla u(\y)) b(\x)v(\x)
\, d\sigma(\y) \, d\sigma(\x).
\end{aligned}
\]
But since $u(\y)v(\y) = 0$ for $|\x - \y| > 20R$, $\|\nabla u\|_{L^\infty(\partial\Omega)} \leq
CR^{-2}$ and $\|v\|_{L^1(\partial\Omega_h)} \leq C$, we can rescale
to obtain
\begin{equation} \label{wbtwo}
\begin{aligned}
& \left| \int_{\partial\Omega_h} \int_{\partial\Omega} (\t(\y) \cdot \nabla
\Gamma_\x(\y)) u(\y) b(\x)v(\x) \,
d\sigma(\y) \, d\sigma(\x) \right| \\
& \leq
CR^{-2} \, \sup_\x \int_{\partial\Omega \cap B_{20R}(\x)} |\Gamma_\x(\y)|
\, d\sigma(\y) \\
& \leq
CR^{-1} \, \int_{-R_1}^{R_1} |\log|y||
\, dy \leq  CR^{-1}.
\end{aligned}
\end{equation}
The same argument as in the end of the proof of Lemma \ref{continuous} shows us
that the left-hand sides of (\ref{wbone}), (\ref{wbthree}) and (\ref{wbtwo}) are all we need to
control in order to prove the weak boundedness of $M_{B_2^t}TM_{B_1}$.

To show that $T(B_1) \in L^\infty(\R,\M)$ we will apply Green's second identity,
\[
\begin{aligned}
\varphi(\x) - \iint_{\Omega_R} {\Gamma}_\x(\y) & L\varphi(\y) \,
d\y = \int_{\partial\Omega_R} \n(\y) \cdot A^t(\y) \nabla {\Gamma}_\x(\y)
\varphi(\y) \, d\sigma(\y) \\
&-
\int_{\partial\Omega_R} {\Gamma}_\x(\y) \n(\y) \cdot A(\y) \nabla \varphi(\y) \,
d\sigma(\y),
\end{aligned}
\]
to the domain $\Omega_R = \Omega \cap B_R(\x)$ and the function $\varphi \equiv
1$. This leads to the equality
\[
\begin{aligned}
1 & = \int_{\partial\Omega_R} \n(\y) \cdot A^t(\y) \nabla {\Gamma}_\x(\y)
\, d\sigma(\y) \\
& = \int_{\partial\Omega \cap B_R(\x)} \n(\y) \cdot A^t(\y) \nabla {\Gamma}_\x(\y)
\, d\sigma(\y) \\
& + \int_{\Omega \cap \partial B_R(\x)} \n(\y) \cdot A^t(\y) \nabla {\Gamma}_\x(\y)
\, d\sigma(\y),
\end{aligned}
\]
but, using the estimate (\ref{103}), we have
\[
\begin{aligned}
\left| \int_{\Omega \cap \partial B_R(\x)} \n(\y) \cdot A^t(\y) \nabla {\Gamma}_\x(\y)
\, d\sigma(\y) \right| \leq C,
\end{aligned}
\]
so, letting $R \to \infty$ in the above, we obtain
\[
\left| \int_{\partial\Omega} \n(\y) \cdot A^t(\y) \nabla {\Gamma}_\x(\y)
\, d\sigma(\y) \right| \leq C.
\]
We also have that, with $\Omega^-_R := \Omega^- \cap B_R(\x)$,
\[
\begin{aligned}
0 & = \int_{\partial\Omega^-_R} \t(\y) \cdot \nabla {\Gamma}_\x(\y)
\, d\sigma(\y) \\
& = \int_{\partial\Omega \cap B_R(\x)} \t(\y) \cdot \nabla {\Gamma}_\x(\y)
\, d\sigma(\y) \\
& + \int_{\Omega^- \cap \partial B_R(\x)} \t(\y) \cdot \nabla {\Gamma}_\x(\y)
\, d\sigma(\y),
\end{aligned}
\]
and in exactly the same fashion as before we can see that
\[
\left| \int_{\partial\Omega} \t(\y) \cdot \nabla {\Gamma}_\x(\y)
\, d\sigma(\y) \right| \leq C.
\]
And so (as in the proof of Lemma \ref{continuous}), we have controlled all the terms necessary to prove
the second statement of the lemma.
\end{proof}

In order to apply Theorem \ref{T(B)David} we need to show that $T^t(B_2) \in
\mbox{BMO}$. At this point the smallness of the Lipschitz constant will come in
(specifically, in Lemma \ref{tildeineq} and in the bound $C_5$). It
will be necessary to first introduce a second singular integral operator
$\widetilde{T}$ defined as
\begin{equation} \label{def2}
\langle G,\widetilde{T}(F)\rangle = \lim_{h \searrow 0} \iint_{\R^2} G(x)^t\widetilde{K}_h(x,y)F(y)
\, dydx,
\end{equation}
with
$\widetilde{K}_h$ defined as follows. For the fundamental solution $\Gamma_\x^t$
of $L^t = \div \, A^t\nabla$ and fixed $\x$ let us define the conjugate
$\widetilde{\Gamma}_\x^t$ of
$\Gamma_\x^t$ to be
\begin{equation} \label{gamtildedef}
\widetilde{\Gamma}_\x^t(\y) = \int_{\gamma(\y_0,\y)} \n(\z) \cdot A(\z) \nabla
\Gamma_\x^t(\z) \,
dl(\z)
\end{equation}
on the complement of the set $\{\y = (y,s) \, | \, s \geq t, \, y = x \}$.
Here $\gamma(\y_0,\y)$ is
a path from a fixed point $\y_0$ to $\y$ parametrised by arc length via the function $t \mapsto (l_1(t), l_2(t))$ and remaining in the
complement of $\{\z = (z,r) \, | \, r \geq t, \, z = x \}$. Also 
$\n(\z) = (l_2'(t),-l_1'(t))$ is the unit normal to $\gamma(\y_0,\y)$ at 
$\z = (l_1(t), l_2(t))$ and $dl$ is arc length. It is easy to see 
$\widetilde{\Gamma}_\x^t(\y)$ solves the system
\begin{equation} \label{system}
A(\y) \nabla \Gamma_\x^t(\y) =
\left(  \begin{array}{cc} 0 & 1 \\ -1 & 0 \end{array} \right)
\nabla \widetilde{\Gamma}_\x^t(\y).
\end{equation}
The function $\y \mapsto \widetilde{\Gamma}_\x^t(\y)$ is well-defined up to a
constant (which depends on the choice of 
$\y_0$). Set
\[
\left(  \begin{array}{c} \widetilde{k}_{1,h}(x,y) \\ \widetilde{k}_{2,h}(x,y)
\end{array} \right) = \nabla_\x \widetilde{\Gamma}^t_{(x,\phi(x)+h)}(y,\phi(y))
\]
and define
\[
\widetilde{K}_h =
\left(  \begin{array}{cc} \widetilde{k}_{1,h} & \widetilde{k}_{2,h} \\
\widetilde{k}_{1,h} & \widetilde{k}_{2,h}
\end{array} \right).
\]
As before we write $\widetilde{K}_0 = \widetilde{K}$.
\begin{lemma} \label{Ttildesingintop}
For each 
$\phi \in \Lambda^\frac{k}{4}(\varepsilon_0)$, the operator $\widetilde{T}$ defined above is a singular integral operator
associated to the Calder\'on-Zygmund kernel $\widetilde{K}$. More precisely,
$\widetilde{K}$ satisfies estimates (\ref{kerest1}), (\ref{kerest2}) and
(\ref{kerest3}) with the constant $C = C_8$ depending only on 
$\lambda$, 
$\Lambda$ and 
$k$, and $\widetilde{T} \colon
B_3\S \to (B_1\S)'$ is a continuous linear operator with representation
(\ref{singrep}).
\end{lemma}

\begin{proof}
The estimates (\ref{kerest1}), (\ref{kerest2}) and (\ref{kerest3})
may be easily obtained from the estimates (\ref{103}), (\ref{104}) and
(\ref{105}) with $\nabla_\x \widetilde{\Gamma}^t_\x$ replacing $\nabla \Gamma_\x$, so it suffices
to show the later. We also observe that, since 
$\phi \in \Lambda^\frac{k}{4}(\varepsilon_0)$, we only need to prove (\ref{103}), (\ref{104}) and (\ref{105}) for $\x = (x,t) \in \{ \z = (z,r) \, | \, |r - s| \leq k/2|z - y| \}$, where $\y = (y,s)$. Since $\x \mapsto
\widetilde{\Gamma}_\x^t(\y)$ is a solution to the non-divergence form equation
$L^tu = 0$ and 
$\y \mapsto \widetilde{\Gamma}_\x^t(\y)$ is a solution to the divergence form equation $\div \, \widetilde{A} \nabla u = 0$ when 
$\x \neq \y$, we may repeat the proof of Lemma \ref{kernelest1} once we have
shown (\ref{103}) for $\nabla_\x \widetilde{\Gamma}_\x^t$. From (\ref{gamtildedef}) we have
\[
\nabla_\x \widetilde{\Gamma}_\x^t(\y) = \int_{\gamma(\y_0,\y)} \nabla_\x (\n(\z) \cdot
A(\z) \nabla \Gamma_\x^t(\z)) \,
dl(\z).
\]
Now, with the a priori smoothness assumption, the integrand here has size
$C/|\x - \z|^2$ for $\z$ away from $\x$, so 
\[
\int_{\gamma(\y_1,\y_2)} \nabla_\x (\n(\z) \cdot
A(\z) \nabla \Gamma_\x^t(\z)) \,
dl(\z) \to 0
\]
as $\min(|\y_1 - \x|, |\y_2 - \x|) \to \infty$, which says we can write
\[
\nabla_\x \widetilde{\Gamma}_\x^t(\y) = \int_{\gamma_\theta(\y)} \nabla_\x (\n(\z) \cdot
A(\z) \nabla \Gamma_\x^t(\z)) \,
dl(\z),
\]
where $\gamma_\theta(\y)$ is a ray from infinity to $\y$ approaching $\y$ at an angle
$\theta$ such that $\arctan (k/2) - \pi \leq \theta \leq -\arctan (k/2)$. Averaging
over $\theta$ we obtain
\[
\begin{aligned}
& \nabla_\x \widetilde{\Gamma}_\x^t(\y) = \fint_{\arctan k -
\pi}^{-\arctan k}\int_{\gamma_\theta(\y)} \nabla_\x (\n(\z) \cdot
A(\z) \nabla \Gamma_\x^t(\z)) \,
dl(\z) d\theta \\
& = C\int_{\Gamma^-_1(\y)} \nabla_\x (\n(\z) \cdot
A(\z) \nabla \Gamma_\x^t(\z)) \,
\frac{d\z}{|\z - \y|} \\
& = C\int_{\Gamma^-_1(\y) \cap B_{|\x - \y|/8}(\y)} \nabla_\x (\n(\z) \cdot
A(\z) \nabla \Gamma_\x^t(\z)) \,
\frac{d\z}{|\z - \y|} \\
& + C\sum_{n = -2}^\infty \int_{\Gamma^-_1(\y) \cap (B_{2^n|\x -
\y|}(\y) \setminus B_{2^{n-1}|\x - \y|}(\y))} \nabla_\x (\n(\z) \cdot
A(\z) \nabla \Gamma_\x^t(\z)) \,
\frac{d\z}{|\z - \y|} \\
& = Q_1 + Q_2,
\end{aligned}
\]
where $\Gamma^-_1(\y) = \{ \z \in \R^2 \,|\, s - r > k |z - y| \}$. Define also $\Gamma^-_2(\y) = \{ \z \in \R^2 \,|\, s - r > 3k/4 |z - y| \}$ and $\Gamma^-_3(\y) = \{ \z \in \R^2 \,|\, s - r > k/2 |z - y| \}$. Observe
that since 
$\x \not\in \Gamma^-_3(\y)$ and 
$\z \in \Gamma^-_1(\y)$ we have $|\x - \z| \geq C(|\x - \y| +
|\z - \y|)$. Thus using
H\"older's inequality, Lemma \ref{cacciopoli} and (\ref{103}), we
can estimate
\[
\begin{aligned}
Q_2 & \leq C \sum_{n = -2}^\infty \int_{\Gamma^-_1(\y) \cap (B_{2^n|\x -
\y|}(\y) \setminus B_{2^{n-1}|\x - \y|}(\y))} |\nabla_\x
\nabla \Gamma_\x^t(\z))| \,
\frac{d\z}{|\z - \y|} \\
& \leq C \sum_{n = -2}^\infty \left( \int_{\Gamma^-_1(\y) \cap (B_{2^n|\x -
\y|}(\y) \setminus B_{2^{n-1}|\x - \y|}(\y))} |\nabla_\x
\nabla \Gamma_\x^t(\z))|^2 \, d\z \right)^\frac{1}{2} \\
& \leq C \sum_{n = -2}^\infty \frac{2^{-n}}{|\x - \y|}\left(
\int_{{\Gamma}^-_2(\y) \cap (B_{2^n|\x -
\y|}(\y) \setminus B_{2^{n-1}|\x - \y|}(\y))} |\nabla_\x \Gamma_\x^t(\z))|^2 \, d\z \right)^\frac{1}{2} \\
& \leq C \sum_{n = -2}^\infty \frac{2^{-n}}{|\x - \y|}\left( \int_{{\Gamma}^-_2(\y) \cap (B_{2^n|\x -
\y|}(\y) \setminus B_{2^{n-1}|\x - \y|}(\y))} \frac{1}{|\x - \z|^2} \, d\z \right)^\frac{1}{2} \\
& \leq C \sum_{n = -2}^\infty \frac{2^{-n}}{|\x - \y|}\left( \int_{{\Gamma}^-_2(\y) \cap (B_{2^n|\x -
\y|}(\y) \setminus B_{2^{n-1}|\x - \y|}(\y))} \frac{1}{(2^{n-1}|\x - \y|)^2} \, d\z \right)^\frac{1}{2} \\
& \leq C \sum_{n = -2}^\infty \frac{2^{-n}}{|\x - \y|} \leq
\frac{C}{|\x - \y|}.
\end{aligned}
\]
Also, using H\"older's inequality, Lemma \ref{reversehoelder} and Lemma
\ref{cacciopoli}, we have, for some $p > 2$,
\[
\begin{aligned}
Q_1 & \leq C \int_{B_{|\x - \y|/8}(\y)} |\nabla_\x
\nabla \Gamma_\x^t(\z))| \,
\frac{d\z}{|\z - \y|} \\
& \leq C \left( \int_{B_{|\x - \y|/8}(\y)} \frac{1}{|\z - \y|^{p'}} \,
d\z \right)^\frac{1}{p'} \left( \int_{B_{|\x - \y|/8}(\y)} |\nabla_\x
\nabla \Gamma_\x^t(\z))|^p \, d\z \right)^\frac{1}{p} \\
& = C |\x - \y|^\frac{2-p'}{p'} \left( \int_{B_{|\x - \y|/8}(\y)} |\nabla_\x
\nabla \Gamma_\x^t(\z))|^p \, d\z \right)^\frac{1}{p} \\
& \leq C \left( \int_{B_{|\x - \y|/4}(\y)} |\nabla_\x
\nabla \Gamma_\x^t(\z))|^2 \, d\z \right)^\frac{1}{2} \\
& \leq \frac{C}{|\x - \y|} \left( \int_{B_{|\x - \y|/2}(\y)} |\nabla_\x
\Gamma_\x^t(\z))|^2 \, d\z \right)^\frac{1}{2} \leq \frac{C}{|\x - \y|}. \\
\end{aligned}
\]
This proves (\ref{103}) for $\nabla_\x \widetilde{\Gamma}_\x^t$ and the other two estimates (\ref{104}) and (\ref{105})
follow as before.

To prove the continuity of $\widetilde{T} \colon
B_3\S \to (B_1\S)'$ we again fix $f, g \in \S(\R, \R)$,
$h \in C^\infty_0(\R,\R)$ with $h$ positive and
equal to one near zero. Let $\delta_0 \colon \Omega^- \to \R$ be the
Dahlberg-Kenig-Stein adapted distance function introduced in \cite{D2}.
Define $u,v \colon \Omega \to \R$ by
\[
u(\x) = f(x)h(\delta_0(x,t)) \hspace{0.2cm} \mbox{and} \hspace{0.2cm} v(\x) =
g(x)h(t - \phi(x)).
\]
Define $\mu \colon \R^2 \to \R$ by $\mu(\x) = (1,\alpha_0)\cdot(x,t)$.
Then $\nabla_\x \mu(\x - \z) = (1, \alpha_0) = \kappa$ for all $\z \in \R^2$.
For $F \colon \Omega^- \to \R$ we define (see \cite{CMS,S}) the quantities
\[
\mathfrak{T}(F)(\q) := \sup_{\Delta_r \ni \q} \left( \frac{1}{|\Delta_r|}
\iint_{B_r(\q)} |F(\y)|^2 \frac{d\y}{\delta(\y)} \right)^\frac{1}{2},
\]
where $\Delta_r = B_r(\q') \cap \partial\Omega$ and $\q' \in \partial\Omega$,
and
\[
\mathfrak{S}(F)(\q) := \left( \iint_{\Gamma^-(\q)} |F(\y)|^2 \,
\frac{d\y}{\delta(\y)^2} \right)^\frac{1}{2}.
\]
We have that
\begin{equation} \label{funcspaces}
\iint_{\Omega^-} |F(\y)G(\y)| \, \frac{d\y}{\delta(\y)} \leq
C\int_{\partial\Omega} \mathfrak{T}(F)(\q) \mathfrak{S}(G)(\q) \, d\sigma(\q).
\end{equation}
It is
well-known that $\|\nabla\delta_0\|_{L^\infty(\Omega,\R)}$ and
$\|\mathfrak{T}(\delta_0\partial^2\delta_0)\|_{L^\infty(\partial\Omega,\R)}$ are bounded in terms of the Lipschitz
constant, where $\partial^2$ denotes any second-order partial derivative, and
$\delta \simeq \delta_0$.

We apply Green's second identity
to the function $\mu(\cdot - \z)u(\cdot)$ and the fundamental
solution ${\Gamma}^t_\x$ in the domain $\Omega^-$. We obtain
\[
\begin{aligned}
& \mu(\x - \z) u(\x) - \iint_{\Omega^-} {\Gamma}^t_\x(\y)
L^t(\mu(\cdot - \z) u(\cdot))(\y)\,
d\y  \\
& = \int_{\partial\Omega} \n(\y) \cdot A(\y) \nabla {\Gamma}^t_\x(\y)
\mu(\y) u(\y) \, d\sigma(\y) \\
& - \int_{\partial\Omega} {\Gamma}^t_\x(\y) \n(\y) \cdot A^t(\y) (u(\y)\kappa +
\mu(\y - \z)\nabla u(\y)) \,
d\sigma(\y).
\end{aligned}
\]
Using (\ref{system}) and integration by parts in the first term on the
right-hand side, we see the above is equal to
\[
\begin{aligned}
& \mu(\x - \z)u(\x) - \iint_{\Omega^-} {\Gamma}^t_\x(\y)
L^t(\mu(\cdot - \z) u(\cdot))(\y))\,
d\y  \\
& = -\int_{\partial\Omega}   \widetilde{\Gamma}^t_\x(\y)
\t(\y) \cdot (u(\y)\kappa + \mu(\y - \z)\nabla u(\y)) \, d\sigma(\y) \\
& - \int_{\partial\Omega} {\Gamma}^t_\x(\y) \n(\y) \cdot A^t(\y)
(u(\y)\kappa +
\mu(\y - \z)\nabla u(\y)) \,
d\sigma(\y).
\end{aligned}
\]
The above integration by parts is the reason for introducing 
$\widetilde{\Gamma}^t_\x$ --- see also (\ref{tilde3}) and (\ref{april1}). Taking the gradient in $\x$, multiplying by $v\b$, where $\b \colon
\partial\Omega_h \to \R^2$ is a bounded
function, setting $\z = \x$ and integrating in $\x$, we find
\begin{equation} \label{tilde1}
\begin{aligned}
& \int_{\partial\Omega_h} u(\x)v(\x)\b(\x) \cdot \kappa \, d\sigma(\x) \\
& - \int_{\partial\Omega_h} v(\x)\b(\x) \cdot \left( \iint_{\Omega^-} \nabla_\x{\Gamma}^t_\x(\y)
L^t(\mu(\cdot - \x) u(\cdot))(\y)\,
d\y \right) d\sigma(\x) \\
& = -\int_{\partial\Omega_h} v(\x)\b(\x) \cdot \int_{\partial\Omega}   \nabla_\x
\widetilde{\Gamma}^t_\x(\y) \times
 \\
& \hspace{2.5cm} \{\t(\y) \cdot (u(\y)\kappa + \mu(\y - \x)\nabla u(\y))\} \, d\sigma(\y) d\sigma(\x)\\
& - \int_{\partial\Omega_h} v(\x)\b(\x) \cdot \int_{\partial\Omega} \nabla_\x
{\Gamma}^t_\x(\y) \times \\
& \hspace{2.5cm} \{\n(\y) \cdot A^t(\y)
(u(\y)\kappa +\mu(\y - \x)\nabla u(\y))\} \,
d\sigma(\y) d\sigma(\x).
\end{aligned}
\end{equation}
The first term on the right-hand side can be expanded:
\begin{equation} \label{tilde4}
\begin{aligned}
& \int_{\partial\Omega_h} v(\x)\b(\x) \cdot \int_{\partial\Omega}   \nabla_\x
\widetilde{\Gamma}^t_\x(\y)u(\y) (\t(\y) \cdot \kappa) \, d\sigma(\y)
d\sigma(\x) \\
& + \int_{\partial\Omega_h} v(\x)\b(\x) \cdot \int_{\partial\Omega}   \nabla_\x
\widetilde{\Gamma}^t_\x(\y)\mu(\y - \x)
\t(\y) \cdot \nabla u(\y))\} \, d\sigma(\y) d\sigma(\x).
\end{aligned}
\end{equation}
The first term in this expression is exactly what we need to control to show the
continuity of $\widetilde{T} \colon B_3\S \to (B_1\S)'$ when $\b$ is chosen to
be either one of the columns of $B_1$. Moreover, the second term on the right-hand side of
(\ref{tilde1}) has already been seen to be controllable by a finite linear
combination of $[f]_N^\beta[g]_N^\beta$, since the proof of Lemma
\ref{continuous} shows that $(T_-)^t \colon B_0\S \to (B_1\S)'$ is continuous for
any bounded $B_0$, where $T_-$ is the operator defined as in (\ref{def}), but
with the limit as $h \nearrow 0$. To control the second term on the left-hand side, we have
\[
|L^t(\mu(\cdot - \x) u(\cdot))| \leq |\mu(\cdot - \x)L^t(u)| + C|\nabla u| 
\]
so
\[
\begin{aligned}
& \iint_{\Omega^-} \nabla_\x{\Gamma}^t_\x(\y)
L^t(\mu(\cdot - \x) u(\cdot))(\y)\,
d\y \\ 
& \leq C\iint_{\Omega^-} |L^t(u)(\y)| \, d\y + C\iint_{\Omega^-} \frac{|\nabla
u(\y)|}{|\x - \y|} \, d\y.
\end{aligned}
\]
The second term is easy to control by $C([f]^1_2 + [f]^0_2)$ and so its
contribution to the second term on the left-hand side of
(\ref{tilde1}) is also controlled. The first term here requires a little more
work. Let us consider just the first term $\partial_{11}u$ of $L^t(u)$, as the same
analysis can be used on the other terms. We have
\[
\begin{aligned}
\partial_{11}u(\x) & = f''(x)h(\delta_0(x,t)) +
f'(x)h'(\delta_0(x,t))\partial_1\delta_0(x,t) \\
& + f(x)h''(\delta_0(x,t))(\partial_1\delta_0(x,t))^2 +
f(x)h'(\delta_0(x,t))(\partial_{11}\delta_0(x,t)),
\end{aligned}
\] 
so the contribution to the double integral from the first three terms is easily
bounded by $C([f]^2_2 + [f]^1_2 + [f]^0_2)$. The last term is controlled using
(\ref{funcspaces}):
\[
\iint_{\Omega^-} |f(y)h'(\delta_0(y,s))(\partial_{11}\delta_0(y,s))| \, d\y
\leq C\int_{\partial\Omega} \mathfrak{T}(f h'(\delta_0))
\mathfrak{S}(\delta_0\partial_{11}\delta_0) \, d\sigma \leq C[f]^0_2,
\]
since
$\|\mathfrak{S}(\delta_0\partial_{11}\delta_0)\|_{L^\infty(\partial\Omega,\R)} \leq C$
and it is clear from the support properties of $h'(\delta)$ that
$\|\mathfrak{S}(fh'(\delta_0))\|_{L^1(\partial\Omega, \R)} \leq C[f]^0_2$. Thus,
the second term on the left-hand side of
(\ref{tilde1}) is also controlled.
The remaining terms in (\ref{tilde1}) and (\ref{tilde4}) may be controlled in
terms of a finite linear combination of $[f]_N^\beta[g]_N^\beta$ directly, as
was done in the proof of Lemma \ref{continuous}. We remark that under our a priori smoothness assumptions the limit (\ref{def2}) can be seen to exist, again as was done in the proof of Lemma \ref{continuous}.
\end{proof}

\begin{remark} \label{awayfromcone} Note that, in fact, we have estimate (\ref{103}) for
$\nabla_\x \widetilde\Gamma^t_\x(\y)$ when $\y \not\in \{ \z \, | \,  r - t >
k/2|z-x| \}$, that is away from a cone with axis being the ray on which
$\widetilde{\Gamma}^t_\x$ is not defined. This will be used in Lemma \ref{Ttildeweakboundedness}.
\end{remark}

\begin{lemma} \label{Ttildeweakboundedness}
For each 
$\phi \in \Lambda^\frac{k}{4}(\varepsilon_0)$, the operator $M_{B_1^t}\widetilde{T}M_{B_3}$ satisfies
the weak boundedness property and
$\widetilde{T}^t(B_1) \in L^\infty(\R,\M)$, with bounds depending only on 
$\lambda$, 
$\Lambda$ and 
$k$.
\end{lemma}

\begin{proof}
We can prove the weak boundedness by using (\ref{tilde1}) with $u$ and $v$
replaced by 
\[
\begin{aligned}
u(\x) & =
R^{-1}f_{ij}(\x/R)h(\delta_0(x,t)/R) \hspace{.2cm} \mbox{and} \\
v(\x) & = R^{-1}g_{ij}(x/R)h((t-\phi(x))/R),
\end{aligned}
\]
with $f_{ij}$, $g_{ij}$ and $h$ as in (\ref{tilde2}),
and $\b$ again being either one of
the columns of $B_1$. The estimate we require is that
\[
\left| \int_{\partial\Omega_h} v(\x)\b(\x) \cdot \int_{\partial\Omega}   \nabla_\x
\widetilde{\Gamma}^t_\x(\y)u(\y) (\t(\y) \cdot \kappa) \, d\sigma(\y)
d\sigma(\x) \right| \leq \frac{C}{R},
\]
which is part of the first term on the right-hand side of (\ref{tilde1}).
The second term on the right-hand side of (\ref{tilde1}) is
controlled in the necessary manner because $M_{B_1^t}(T_-)^tM_{B_0}$ satisfies the
weak boundedness property for any bounded $B_0$ (as the proof of Lemma
\ref{weakboundedness} shows). The second term on the left-hand side of
(\ref{tilde1}) is controlled in the same manner as in the proof of Lemma
\ref{Ttildesingintop}, and the remaining terms can be
estimated directly, leaving us with the required estimate.

To show that $\widetilde{T}^t(B_1) \in L^\infty(\R,\M)$ we will apply Green's second identity,
\[
\begin{aligned}
- \iint_{\Omega_{h,R}} \widetilde{\Gamma}_\x(\y) & L\varphi(\x) \,
d\x = \int_{\partial\Omega_{h,R}} \n(\x) \cdot A^t(\x) \nabla_\x \widetilde{\Gamma}_\x(\y)
\varphi(\x) \, d\sigma(\x) \\
& -
\int_{\partial\Omega_{h,R}} \widetilde{\Gamma}_\x(\y) \n(\x) \cdot A(\x) \nabla \varphi(\x) \,
d\sigma(\x),
\end{aligned}
\]
to the domain $\Omega_{h,R} := \{\x = (x,t) \, | \, t > \phi(x) + h \} \cap
B_R(\y)$ and the functions $\x \mapsto \widetilde{\Gamma}_\x^t(\y)$ and $\varphi \equiv 1$. This leads to the equality
\[
\begin{aligned}
0 & = \int_{\partial\Omega_{h,R}} \n(\x) \cdot A^t(\x) \nabla_\x \widetilde{\Gamma}_\x(\y)
\, d\sigma(\x) \\
& = \int_{\partial\Omega_h \cap B_R(\y)} \n(\x) \cdot A^t(\x) \nabla_\x \widetilde{\Gamma}_\x(\y)
\, d\sigma(\x) \\
& + \int_{\Omega_h \cap \partial B_R(\y)} \n(\x) \cdot A^t(\x) \nabla_\x \widetilde{\Gamma}_\x(\y)
\, d\sigma(\x),
\end{aligned}
\]
but, by Remark \ref{awayfromcone}, we have
\[
\begin{aligned}
\left| \int_{\Omega_h \cap \partial B_R(\x)} \n(\x) \cdot A^t(\x) \nabla_\x \widetilde{\Gamma}_\x(\y)
\, d\sigma(\x) \right| \leq C,
\end{aligned}
\]
so, letting $R \to \infty$ in the above, we obtain
\[
\left| \int_{\partial\Omega} \n(\x) \cdot A^t(\x) \nabla_\x \widetilde{\Gamma}_\x(\y)
\, d\sigma(\x) \right| \leq C.
\]
The bound
\[
\left| \int_{\partial\Omega_h} \t(\x) \cdot \nabla_\x \widetilde{\Gamma}_\x^t(\y)
\, d\sigma(\x) \right| \leq C
\]
also follows via the method used to
prove Lemma \ref{weakboundedness} on $\Omega_{h,R}$. These two estimates bound all the components of $\widetilde{T}^t(B_1)$.
\end{proof}

\begin{lemma} \label{tildeineq}
For each 
$\phi \in \Lambda^\frac{k}{4}(\varepsilon_0)$, we have the inequalities $\|\widetilde{T}(B_3)\|_{\mbox{\tiny BMO}} \leq C +
C\|T^t(B_2)\|_{\mbox{\tiny BMO}}$ and $\|T^t(B_2)\|_{\mbox{\tiny BMO}} \leq C +
C\varepsilon_0\|\widetilde{T}(B_3)\|_{\mbox{\tiny BMO}}$, with the constants 
$C$ depending only on 
$\lambda$, 
$\Lambda$ and 
$k$.
\end{lemma}

\begin{proof}
To prove this lemma we will use the $H^1$-BMO duality and the smooth atomic decomposition of $H^1$ (for a discussion of this see, for example, \cite{G}). For each $H^1$ atom $A_0 \in B_1 \S$, with support in a ball $B_R(x_0)$ of radius
$R$, we choose a $\eta \colon \R \to \R$ with support in $B_{4R}(x_0)$ and identically
equal to one on $B_{2R}(x_0)$ such that $|\eta'|
\leq C/R$ and $|\eta''|
\leq C/R^2$. Recall that with our a priori smoothness assumptions such an atom 
$A_0$ is in fact an arbitrary smooth atom. To prove the first inequality we use definition (\ref{defBMO}) to
compute $\langle A_0, \widetilde{T}(B_3) \rangle$.
The second term on the right-hand side of (\ref{defBMO}) is
controlled by a multiple of $\|B_3\|_\infty$, and so we only
need to estimate $\langle A_0, \widetilde{T}(\eta B_3) \rangle$. To that end, we define $\varphi \colon \R^2 \to
\R$ by
\begin{equation} \label{varphi}
\varphi(\x) = \eta(x)\eta(\delta_0(x,t)),
\end{equation}
where we recall $\x = (x,t)$ and $\delta_0$ is the adpated distance function,
and repeat the calculation that yielded (\ref{tilde1}) with
$\mu(\cdot - \z)\varphi(\cdot)$ and $\tilde{a}$, where $\tilde{a}(x,\phi(x) + h) =
a(x)/(1+\phi'(x))^\frac{1}{2}$ and $a$ is any matrix element of $A_0$. This
gives
\begin{equation} \label{tilde3}
\begin{aligned}
& \int_{\partial\Omega_h} \tilde{a}(\x)\varphi(\x) \kappa \, d\sigma(\x) \\
& + \int_{\partial\Omega_h}  \tilde{a}(\x) \left( \iint_{\Omega^-}
\nabla_\x{\Gamma}^t_\x(\y)
L^t(\mu(\cdot - \x) \varphi(\cdot))(\y)\,
d\y \right) d\sigma(\x) \\
& = \int_{\partial\Omega_h} \tilde{a}(\x) \int_{\partial\Omega}   \nabla_\x
\widetilde{\Gamma}^t_\x(\y) \times
 \\
& \hspace{2.5cm} \{\t(\y) \cdot (\varphi(\y)\kappa + \mu(\y - \x)\nabla \varphi(\y))\} \, d\sigma(\y) d\sigma(\x)\\
& + \int_{\partial\Omega_h} \tilde{a}(\x) \int_{\partial\Omega} \nabla_\x
{\Gamma}^t_\x(\y) \times \\
& \hspace{2.5cm} \{\n(\y) \cdot A^t(\y)
(\varphi(\y)\kappa +\mu(\y - \x)\nabla \varphi(\y))\} \,
d\sigma(\y) d\sigma(\x).
\end{aligned}
\end{equation}

The first term on the left-hand side of (\ref{tilde3}) is easy to control.
Indeed
\[
\begin{aligned}
& \int_{\partial\Omega_h} \tilde{a}(\x)\varphi(\x) \kappa \, d\sigma(\x) \\
& \leq C\|a\|_{L^2(\R)} \|\varphi\kappa\|_{L^2(\partial\Omega_h)} \leq
CR^\frac{1}{2}/R^\frac{1}{2}.
\end{aligned}
\]
To control the second term on the left-hand side of (\ref{tilde3}) we expand
\[
\begin{aligned}
L^t(\mu\varphi) & = \varphi L^t\mu + \mu L^t\varphi + 2\partial_1\varphi
\partial_1 \mu + b(\partial_1\varphi
\partial_2 \mu + \partial_2\varphi
\partial_1 \mu) + 2c\partial_2\varphi
\partial_2 \mu \\
\end{aligned}
\]
The first term here is zero and the other terms except for
$\mu L^t\varphi$ are of size 
$C/R$. Now
\[
\begin{aligned}
\partial_{11}\varphi(x,t)  & = \eta''(x) \eta(\delta_0(x,t)) + 2\eta'(x) \eta'(\delta_0(x,t)) \partial_1\delta_0(x,t) \\ 
& + \eta(x) \eta''(\delta_0(x,t)) (\partial_1\delta_0(x,t))^2 + \eta(x) \eta'(\delta_0(x,t)) \partial_{11}\delta_0(x,t),
\end{aligned}
\]
so each term here is of size 
$C/R^2$ except for the term where two derivatives fall on 
$\delta_0$. Since the same is true for the other second-order derivatives of 
$\varphi$ we can say $\mu L^t\varphi$ is a sum of terms of size 
$C/R$ plus terms involving the second-order derivatives of 
$\delta_0$. We conclude that
\[ \label{BMO4}
\begin{aligned}
L^t(\mu\varphi) & =  \eta \, \eta'(\delta_0)\sum_{i,j = 1,2} c_{ij}\partial_{ij}\delta_0 + Q,
\end{aligned}
\]
where $Q$ is a sum of terms of size $C/R$ and $c_{ij}$ ($i,j = 1,2$) are
functions of size $C|\y - \x|$.
So for $0 < h < \min\{1,R/2\}$,
\[
\left| \iint_{\Omega^- \cap \{R < \delta(\y) < 4R\}} \nabla_\x{\Gamma}^t_\x(\y)Q(\y)
\, dY \right| \leq \frac{CR^2}{RR} \leq C.
\]
Moreover, when $0 < \delta(\y) < R$ on $\Omega^-$ the supports of $\tilde{a}$ and $\nabla \varphi$
are disjoint and separated by a distance $CR$, so it follows that the same is
true of $\tilde{a}$ and
$Q$, thus
$|\nabla_\x\Gamma^t_\x(\y)| \leq C/R$, so once again
\[
\left| \iint_{\Omega^- \cap \{0 < \delta(\y) < R\}} \nabla_\x{\Gamma}^t_\x(\y)Q(\y)
\, dY \right| \leq \frac{CR^2}{RR} \leq C.
\]
Since the integrand in the integral over $\Omega^-$ in (\ref{tilde3}) is
supported on $4B$, the contribution from $Q$ to the second term on the left-hand
side of (\ref{tilde3}) is now controlled. The contribution from $\eta \,
\eta'(\delta_0)\sum c_{ij}\partial_{ij}\delta_0$ is again controlled
using (\ref{funcspaces}) as in the proof of Lemma \ref{Ttildesingintop}.

We will now consider the first term on the right-hand side of (\ref{tilde3}):
\[
\begin{aligned}
& \int_{\partial\Omega_h} \tilde{a}(\x) \int_{\partial\Omega}   \nabla_\x
\widetilde{\Gamma}^t_\x(\y) \t(\y) \cdot \kappa \varphi(\y) d\sigma(\y) d\sigma(\x) \\
& + 
\int_{\partial\Omega_h} \tilde{a}(\x) \int_{\partial\Omega}   \nabla_\x
\widetilde{\Gamma}^t_\x(\y) \t(\y) \cdot \nabla
\varphi(\y) \mu(\y - \x) \, d\sigma(\y) d\sigma(\x).
\end{aligned}
\]
The first term is a component of $\langle A_0, \widetilde{T}(\eta B_3) \rangle$
and the second is easily seen to be bounded by a constant using (\ref{103}). In
exactly the same manner, the second term on the right-hand side of
(\ref{tilde3}) is bounded by $|\langle A_0, T^t(\eta B_0) \rangle| + C$, for some
bounded $B_0$. All
this yields the estimate
\[
|\langle A_0, \widetilde{T}(\eta B_3) \rangle| \leq C + C |\langle A_0,
T^t(\eta B_0)\rangle|,
\]
from which duality, standard Calder\'on-Zygmund theory, Theorem
\ref{T(B)David} and Lemma \ref{weakboundedness} allow us to conclude
\[
\|\widetilde{T}(\eta B_3)\|_{\mbox{\tiny BMO}} \leq C + C\|T^t(B_2)\|_{\mbox{\tiny BMO}}.
\]
The first inequality of the lemma follows from this.

The second inequality follows by the same analysis, replacing $\mu(\cdot - \z)$
by $\xi(\cdot - \z)$ in the calculation that gave (\ref{tilde3}), where $\xi(\x) =
(-\alpha_0,1)\cdot(x,t)$. The first term on the
right-hand side of the equation which replaces (\ref{tilde3}) is
\begin{equation} \label{april1}
\begin{aligned}
& \int_{\partial\Omega_h} \tilde{a}(\x) \int_{\partial\Omega}   \nabla_\x
\widetilde{\Gamma}^t_\x(\y)(\t(\y) \cdot \kappa^\perp) \varphi(\y) \, d\sigma(\y)
d\sigma(\x) \\
& + 
\int_{\partial\Omega_h} \tilde{a}(\x) \int_{\partial\Omega}   \nabla_\x
\widetilde{\Gamma}^t_\x(\y)(\t(\y) \cdot\nabla
\varphi(\y)) \xi(\y - \x) \, d\sigma(\y) d\sigma(\x).
\end{aligned}
\end{equation}
We have that $(1+(\phi')^2)^\frac{1}{2} \t \cdot \kappa^\perp = \phi' - \alpha_0$ which
is bounded in absolute value by $\varepsilon_0$, so the first term in
(\ref{april1}) is a
component of $\langle A_0, \widetilde{T}(B_{\varepsilon_0})\rangle$, where
$B_{\varepsilon_0}$ is a function bounded by $C\varepsilon_0$. The second term
in (\ref{april1}) is bounded by a constant as before. The second term on the
right-hand side of the equation which replaces (\ref{tilde3}) is
\[
\begin{aligned}
& \int_{\partial\Omega_h} \tilde{a}(\x) \int_{\partial\Omega} \nabla_\x
{\Gamma}^t_\x(\y)(\n(\y) \cdot A^t(\y) \kappa^\perp)
\varphi(\y) \,
d\sigma(\y) d\sigma(\x) \\
& + 
\int_{\partial\Omega_h} \tilde{a}(\x) \int_{\partial\Omega} \nabla_\x
{\Gamma}^t_\x(\y) (\n(\y) \cdot A^t(\y)\nabla \varphi(\y))\xi(\y - \x) \,
d\sigma(\y) d\sigma(\x).
\end{aligned}
\]
The first term here is a component of $\langle A_0, T^t(\eta B_2)\rangle$ and the
second is again bounded by a constant, so we have from Calder\'on-Zygmund theory, Theorem \ref{T(B)David} and Lemma \ref{Ttildeweakboundedness} that
\[
\begin{aligned}
\|T^t(\eta B_2)\|_{\mbox{\tiny BMO}} & \leq C + (C +
C\|\widetilde{T}(B_3)\|_{\mbox{\tiny BMO}})\|B_{\varepsilon_0}\|_{L^\infty(\R)} \\
& \leq C + C\varepsilon_0\|\widetilde{T}(B_3)\|_{\mbox{\tiny BMO}},
\end{aligned}
\]
which leads to the second inequality.
\end{proof}

It is now straight forward to complete the proof of Theorem \ref{T(B)}. Lemma
\ref{tildeineq} allows us to conclude
\[
\|T^t(B_2)\|_{\mbox{\tiny BMO}} \leq C +
C\varepsilon_0\|\widetilde{T}(B_3)\|_{\mbox{\tiny BMO}} \leq C + C\varepsilon_0(1 +
\|T^t(B_2)\|_{\mbox{\tiny BMO}})
\]
so we can choose $\varepsilon_0$ sufficiently small, depending only on $k$ and the
ellipticity constants, to hide the BMO-norm on the left-hand side and conclude
\begin{equation} \label{end1}
\|T^t(B_2)\|_{\mbox{\tiny BMO}} \leq C.
\end{equation}
Note that we must check that our a priori assumptions imply that 
$\|T^t(B_2)\|_{\mbox{\tiny BMO}} < \infty$ in order to justify the last step. We can prove this quite quickly in our situation. First, observe that, via Theorem \ref{T(B)David} and standard Calder\'on-Zygmund theory, it sufficies to show $\|T^t(B_4)\|_{\mbox{\tiny BMO}} < \infty$ for some other bounded 
$B_4$ with bounded inverse. Also we have that 
$y \mapsto \phi(y) - \alpha_0y \in C^\infty_0(\R)$ for some 
$\alpha_0 \in \R$, so 
$\|\phi(y) - \alpha_0y\|_{L^\infty(\R)} = M < \infty$. Choose 
$\eta \in C^\infty_0(\R)$ so that 
$\eta(s) = 1$ for 
$|s| \leq 2M$. Setting $d_0(y,s) = s - \alpha_0y - \eta(s - \alpha_0y)(\phi(y) - \alpha_0y)$ and, applying Green's second identity to $\Gamma^t_\x$ and $d_0$, we obtain
\[
d_0(\x) - \iint_\Omega \Gamma^t_\x L^t(d_0)(\y) \, d\y = \int_{\partial\Omega} \Gamma^t_\x(\y) \n(\y) \cdot A^t(\y) \nabla d_0(\y) \, d\sigma(\y),
\]
since 
$d_0 \equiv 0$ on 
$\partial\Omega$. Taking the gradient in 
$\x$ we have
\[
\nabla d_0(\x) - \iint_\Omega \nabla_\x\Gamma^t_\x L^t(d_0)(\y) \, d\y = \int_{\partial\Omega} \nabla_\x \Gamma^t_\x(\y) (\n(\y) \cdot A^t(\y) \nabla d_0(\y)) \, d\sigma(\y).
\]
It is easy to check using (\ref{103}) that both terms on the left are bounded, thus so is the right-hand side. But, if we choose 
$B_4 \colon \R \to \M$ to be the diagonal matrix with diagonal entries both being 
$\n \cdot A^t \nabla d_0$, then control of the right-hand side above controls the columns of $T^t(B_4)$, and 
$B_4$ is indeed bounded with a bounded inverse. Thus $\|T^t(B_4)\|_{L^\infty(\R)} < \infty$ and so $\|T^t(B_2)\|_{\mbox{\tiny BMO}} < \infty$, as required.

We can then use (\ref{end1}) and Lemma \ref{weakboundedness} to apply Theorem \ref{T(B)David} and obtain the
desired bound on the $L^2$ operator norm of $T$ (and
$\widetilde{T}$ via Lemmata \ref{Ttildesingintop}, \ref{Ttildeweakboundedness} and \ref{tildeineq}), completing the proof of Theorem \ref{T(B)}.

\section{Boundedness of the layer potentials on boundaries with arbitrary
Lipschitz constants}

\label{david} 

The aim of this section is to remove the necessity for $\varepsilon_0$ in
Theorem \ref{T(B)} to be small. We formulate this as the theorem below.  It is
proved by applying the build-up scheme of David \cite{D}.

\begin{theorem} \label{genT(B)}
The conclusion of Theorem \ref{T(B)} holds with $\alpha_0 = 0$ and
$\varepsilon_0 = k/8$, that is, the conclusion holds for an arbitrary
Lipschitz function $\phi$.
\end{theorem}

Theorem \ref{genT(B)} will proved using the following (\cite[Prop 10]{D} and \cite[p110]{J}, respectively).

\begin{theorem} \label{sunrise}
Let $\varepsilon_0 > 0$, $\phi \in \Lambda^k(\varepsilon_0)$ be
such that $\|\phi' - \alpha_0\|_{L^\infty(\R)} \leq \varepsilon_0$, and $I \subset \R$ be an
interval. Then there exists a compact subset $E \subset I$ and a function $\psi \in
\Lambda^{k + \frac{\varepsilon_0}{10}}(9\varepsilon_0/10)$ such that
\[
\begin{aligned}
& |E| \geq \frac{1}{3(1 + (k + \varepsilon_0)^2)^\frac{1}{2}}|I|, \\
& \phi(x) = \psi(x) \hspace{.2cm} \mbox{for all} \hspace{.2cm} x \in E, 
\hspace{.2cm} \mbox{and either} \\
& -\frac{4}{5}\varepsilon_0 \leq \psi'(x) - \alpha_0 \leq \varepsilon_0
\hspace{.2cm} \mbox{or} \\
& -\varepsilon_0 \leq \psi'(x) - \alpha_0 \leq \frac{4}{5}\varepsilon_0
\hspace{.2cm} \mbox{almost
everywhere.}
\end{aligned}
\]
\end{theorem}

\begin{theorem} \label{build-up}
Let $K \colon \R^2 \to \M$ be a Calder\'on-Zygmund kernel. Suppose that there exist constants 
$\theta \in (0,1]$ and
$C_{9} > 0$ such that the constant $C$ appearing
in (\ref{kerest1}), (\ref{kerest2}) and (\ref{kerest3}) is no more than $C_{9}$ and, for all intervals $I$, there exists a compact subset $E \subset I$ and a Calder\'on-Zygmund kernel 
$K_I \colon \R^2 \to \M$ have the following properties:
\[
\begin{aligned}
& |E| \geq \theta |I|; \\
& \mbox{for all $x, y \in E$ we have} \hspace{.2cm}  
K_I(x,y) = K(x,y); \\
& \mbox{and} \hspace{.2cm} \|T_I^*\|_{L^2(\R) \to L^2(\R)} \leq C_{9},
\end{aligned}
\]
where $T_I^*$ is the maximal singular integral operator associated to 
$K_I$. 
Then the maximal singular integral operator
$T^*$ associated to 
$K$ is bounded on 
$L^2(\R)$ with $\|T^*\|_{L^2(\R) \to L^2(\R)} \leq C(\theta) C_{9}$.
\end{theorem}

\begin{proof3}
It suffices to show that Theorem \ref{T(B)} holds for $\phi \in
\Lambda^\frac{k}{8}(k/8)$. To this end, with 
$\varepsilon_0$ as in Theorem \ref{T(B)}, pick $m$ so large that $(9/10)^{m} k/8 < \varepsilon_0$. Then Theorem \ref{T(B)} holds for 
$\phi \in \Lambda^{a_mk}((9/10)^{m} k/8)$, where 
$a_m = 1/4$.

We now claim Theorem \ref{T(B)} holds for 
$\phi \in \Lambda^{a_{m-1} k}((9/10)^{m-1} k/8)$, where 
$a_{m-1} = 1/4 -  1/80(9/10)^{m-1}$. To see this fix $\phi \in \Lambda^{a_{m-1} k}((9/10)^{m-1} k/8)$, then for any interval 
$I$ apply Theorem \ref{sunrise} to obtain a compact subset 
$E \subset I$ and a Lipschitz function 
$\phi_I \in \Lambda^{a_mk}((9/10)^{m} k/8)$. Denote by $K_{I}$, the
kernel obtain via (\ref{kernel}) as $K$ was but with $\phi$ replaced with
$\phi_I$. Then $K_{I}$ satisfies the conditions of Theorem \ref{build-up} with $\theta = 1/(3(1 +
k^2)^\frac{1}{2})$. Indeed, Corollary \ref{kernelest2} tells us that
$K_{I}$ is a Calder\'on-Zygmund
kernel, Theorem \ref{T(B)} tells us that
$T_{I}^*$ is bounded on
$L^2(\R)$ and the remaining properties follow from Theorem \ref{sunrise}.

We can now repeat this argument to show that Theorem \ref{T(B)} holds for $\phi \in \Lambda^{a_{m-2} k}((9/10)^{m-2} k/8)$, where 
$a_{m-2} = 1/4 -  1/80((9/10)^{m-1} + (9/10)^{m-2})$. Continuing in this way, after 
$m$ steps we see that Theorem \ref{T(B)} holds for $\phi \in \Lambda^{a_{0} k}(k/8)$, where
\[
a_0 = \frac{1}{4} -  \frac{1}{80}\sum_{j=0}^{m-1}\left(\frac{9}{10}\right)^{j} > \frac{1}{4} -  \frac{1}{80}\sum_{j=0}^{\infty}\left(\frac{9}{10}\right)^{j} = \frac{1}{8}.
\]
Thus Theorem \ref{genT(B)} is proved, and with it our main result, Theorem
\ref{main}.
\end{proof3}

\appendix \label{appendix}

\section*{Appendix: An Example}

Here we will show that given any $p > 1$ there exist operators for which $(R)_p$ and $(N)_p$ do not hold.

Recall the operators $L_h = \div \, A\nabla \cdot$, for $h > 0$, appearing in \cite[Thm (3.2.1)]{KKPT} where
\[
A = \left(
\begin{array}{ll} 1 & m(x) \\ -m(x) & 1 \end{array}
\right)
\]
and
\[
m(x) = \left\{ \begin{array}{ll} h, & \hspace{0.2cm} x \geq 0 \\
-h, & \hspace{0.2cm} x<0.
\end{array} \right.
\]
It is shown there that $u$ solves $L_hu = 0$ in $\R^2_+$ if and only if $u$ is harmonic in the quarter planes $\{\x = (x,t) \, | \, x>0, t>0 \}$ and $\{\x = (x,t) \, | \, x<0, t>0 \}$, smooth up to the boundary except at $(0,0)$, continuous at $(0,0)$ and satisfies the transmission condition
\[
u_x^- - u_x^+ - 2hu_t = 0 \, \, \mbox{on} \, \, \{(x,t) \, | \, x=0\}.
\]
As a result they show that the function $w \colon \R^2_+ \to \R$ given by
\[
w(x,t) = \left\{ \begin{array}{ll} \Im ((x + it)^a), & \hspace{0.2cm} x \geq
0 \\ 
\Im ((-x + it)^a), & \hspace{0.2cm} x<0
\end{array} \right.
\]
satisfies $L_hw = 0$ if and only if $h = \tan b\pi/2$, where $b = 1 - a$. We can readily check that
\[
w = 0, \, \partial_xw = 0, \hspace{.2cm} \mbox{and} \hspace{.2cm} \partial_tw = a|x|^{-b} \hspace{.2cm}
\mbox{on} \hspace{.2cm} \partial\R^2_+
\]
Now solve (\ref{dirichlet}) with data $f_0 \in C^\infty_0(\partial\R^2_+,\R)$ such that $f_0 = 0$ for $|x| < 1$ and $|x| > 2$, $f_0 = 1$ for $9/8 < |x| < 15/8$, and $f_0 \geq 0$ to obtain $u \in \widetilde{W}^{1,2}(\R^2_+)$ via Lemma \ref{dirichletexist}. An application of the comparison principle \cite[Lem 1.3.7]{K} shows that on $\{(x,t) \, | \, |x|<1/2, t=0\}$ we have
\[
\partial_tu \simeq \partial_tw,
\]
thus $\partial_tu \simeq |x|^{-b}$ on the same set. Now, if $bp > 1$, we have $\|\partial_tu\|_{L^p(\partial\R^2_+)} = \infty$ and the regularity of the coefficients ensures that $\nabla u$ converges everywhere on the boundary except perhaps at $(0,0)$, so $\|\widetilde{N}(\nabla u)\|_{L^p(\partial\R^2_+)} = \infty$. However, $\|\partial_xu\|_{L^p(\partial\R^2_+)} < \infty$, so given any $p>1$ we can certainly find $b$ and $h$ so that $(R)_p$ does not hold for $L_h$ in $\R^2_+$.

We can now show that $(N)_p$ cannot hold for the conjugate operator, which has coefficient matrix $\widetilde{A} = A^t/\det(A)$. Let $\widetilde{u}$ be the conjugate of $u$ defined by (\ref{conj2}). Since $(u,\widetilde{u})$ satisfies (\ref{conj}), $\|\nabla u\|_{L^2(\Omega)} < \infty$ and again the regularity of the coefficients ensures that $\widetilde{u}$ is the unique solution to (\ref{neumann}) with data $\n \cdot \widetilde{A} \nabla \widetilde{u}$. Since the conormal derivative of $u$ becomes the tangential derivative of $\widetilde{u}$ and vice versa, $\|\n \cdot \widetilde{A} \nabla \widetilde{u}\|_{L^p(\partial\R^2_+)} < \infty$ but $\|\partial_x\widetilde{u}\|_{L^p(\partial\R^2_+)} = \infty$, therefore, as before, we see $(N)_p$ cannot hold for this operator.

\end{document}